\documentclass[11pt]{amsart}
\usepackage{amssymb}
\usepackage{graphicx}
\usepackage{slashed}
\usepackage[colorlinks]{hyperref} 
\usepackage[all,cmtip]{xy}

\newtheorem{thm}{Theorem}[subsubsection]
\newtheorem{cor}[thm]{Corollary}
\newtheorem{thm-def}[thm]{Theorem-Definition}
\newtheorem{lem}[thm]{Lemma}

\newtheorem{defn}[thm]{Definition}
\theoremstyle{remark}
\newtheorem{rem}[thm]{Remark}
\theoremstyle{remark}

\numberwithin{equation}{subsection}

\newcommand{\vacA}{\mathbf{1}}

\newcommand{\CC}{\mathbb C}
\newcommand{\HH}{\mathbb H}
\newcommand{\QQ}{\mathbb Q}
\newcommand{\ZZ}{\mathbb Z}

\newcommand{\PP}{\mathbb P}
\newcommand{\SSS}{\mathbb S}

\newcommand{\calD}{{\mathcal{D}}}

\newcommand{\cA}{\mathcal{A}}

\newcommand{\cD}{{\mathcal{D}}}
\newcommand{\cE}{{\mathcal{E}}}

\newcommand{\cH}{{\mathcal{H}}}

\newcommand{\cN}{\mathcal{N}}
\newcommand{\cO}{{\mathcal{O}}}

\newcommand{\cT}{{\mathcal{T}}}
\newcommand{\cV}{{\mathcal{V}}}
\newcommand{\cW}{{\mathcal{W}}}

\newcommand{\cZ}{{\mathcal{Z}}}

\newcommand{\ops}[1]{{}_{({#1})} }

\newcommand{\fg}{\frak{g}}

\newcommand{\fm}{\frak{m}}
\newcommand{\fn}{\frak{n}}

\newcommand{\fz}{\frak{z}}

\newcommand{\fU}{{\frak{U}}}

\begin{document}

\title{localization of affine w-algebras}

\author{ T.~Arakawa\and T.~Kuwabara\and F.~Malikov}
\thanks{T. Arakawa is partially   supported 
by the JSPS Grant-in-Aid  for Scientific Research (B)
No.\ 20340007
and the JSPS Grant-in-Aid for challenging Exploratory Research
No.\ 23654006}
\thanks{T.~Kuwabara was partially supported by Basic Science Research Program through the National Research Foundation of Korea (NRF) grant funded by the Korea
government (MEST)(2011-0027952). }

\thanks{F.~Malikov is partially supported by an NSF grant}

\maketitle

\noindent

\begin{abstract}
We introduce the notion of an asymptotic algebra of chiral differential operators. We then construct, via a chiral Hamiltonian reduction, one such
algebra over a resolution of the intersection of the Slodowy slice with the nilpotent cone. We compute the space of global sections of this algebra
thereby proving a localization theorem for affine W-algebras at the critical level.
\end{abstract}

\section{introduction}
\label{introduction}
\subsection{ }
\label{intro-1}
If $X$ is a scheme, $\cA$ a sheaf of associative unital algebras on $X$, $A$ is an associative unital algebra s.t. $A=\Gamma(X,\cA)$, then
one says that $\cA$ is a {\em localization} of $A$. 

If $A$ is commutative, then a certain localization of $A$ is a cornerstone of the very definition of the affine scheme $\text{Spec}A$.

An important noncommutative example was introduced simultaneously in 2 classic works by Beilinson-Bernstein  \cite{BB1,BB2} and Brylinski-Kashiwara \cite{BrK}.
In this case $X$ is a flag variety, and $\cD_X$, the algebra of differential operators on $X$, is a localization of the `centerless'  universal enveloping algebra, $U(\fg)_0$, of the
corresponding simple complex Lie algebra $\fg$. Introduced as a tool to prove the Kazhdan-Lusztig conjecture, this construction has transcended its original
representation-theoretic purpose and served as a template mathematically to define  conformal field theory, see the unpublished but influential \cite{BFM} and the more recent
\cite{BD}.

Relatively recently more examples of noncommutative localization have been considered in \cite{KR}, \cite{Losev} and then \cite{DK}.  To formulate (part of) the Dodd-Kremnizer result
consider the nilpotent cone $\cN\subset\fg$, the Slodowy slice $\SSS\subset\fg$, and the Springer resolution 
\begin{equation}
\label{spr-res-first-app}
\pi: T^*X\rightarrow\cN.
\end{equation}
 We then obtain $S:=\SSS\cap\cN$ and
$\tilde{S}\subset T^*X$, the preimage of $S$ under (\ref{spr-res-first-app}). Dodd and Kremnizer assert that a certain (constructed via Hamiltonian reduction)
 algebra of
asymptotic differential operators over $\tilde{S}$, $\cD_{\tilde{S},\sqrt{\hbar}}$,  is a localization of the asymptotic version of the centerless
finite W-algebra; the result for the genuine W-algebra looks also attractive: $\cW^{fin}_0=\Gamma(\tilde{S},\cD_{\tilde{S},\sqrt{\hbar}})^{\CC^*}$.

We have used the word `asymptotic' and hinted, via the superscript $\CC^*$, at the $\CC^*$-action; these reflect two important constructions invented by Kashiwara and
Rouquier in the above cited \cite{KR} for the purposes of localizing the Cherednik algebra.  To put these constructions in a suitable context let us return to the Beilinson/Bernstein-Brylinski/Kashiwara situation.  Kostant's theorem asserts  that although (\ref{spr-res-first-app}) is only a birational isomorphism, it defines an isomorphism of the rings of regular functions
\begin{equation}
\label{konst-thm-first-app}
\CC[\cN]\stackrel{\sim}{\longrightarrow}\Gamma(T^*X,\cO_{T^*X}).
\end{equation}
This result is the key to the localization isomorphism $U(\fg)_0\stackrel{\sim}{\longrightarrow}\Gamma(X,\cD_X)$. One can say that the latter is a quantization of the
former because, as everybody knows, the algebra of differential operators is a quantization of the algebra of functions on the cotangent bundle.  Everybody knows this, but
the assertion is not quite correct.  A quantization of $\cO_{T^*X}$ is to be understood as a sheaf over $T^*X$, and $\cD_X$, being noncommutative, resists localization along
the fibers of the projection $T^*X\rightarrow X$. Therefore $\cD_X$ is actually a quantization of the push-forward of $\cO_{T^*X}$ on $X$. In  case of the cotangent bundle
this subtlety can be ignored, but the problem is that although $\tilde{S}\subset T^*X$ is symplectic it is not the cotangent bundle of anything.  The way around was suggested
in \cite{KR}. First $\cD_X$ is to be replaced by its deformation quantization version. This means introducing the `Planck constant' $\hbar$ so as to deform the
local defining relations
as follows: $[\xi,\eta]_\hbar=\hbar[\xi,\eta]$, $[\xi,f]_\hbar=\hbar\xi(f)$, $\xi,\eta\in\cT_X$, $f\in\cO_X$;  then complete by the positive powers of $\hbar$ and allow the division by $\hbar$. The result is a
topological algebra of what we called  asymptotic differential operators. It localizes well,\footnote{this phenomenon is best understood from the point of view developed 
by Kapranov in \cite{Kapr}, an important source of inspiration for us} but it is too huge to reproduce a reasonable associative algebra. To the rescue
comes the $\CC^*$-equivariant structure discovered  in \cite{KR} in the case of the Cherednik algebra. Dodd and Kremnizer use this technology in the W-algebra case, where
the $\CC^*$-equivariant structure is a direct manifestation of the well-known Kazhdan filtration. Note an analogue of  (\ref{konst-thm-first-app}):
\begin{equation}
\label{konst-thm-first-app-for-slodow-intro}
\CC[S]\stackrel{\sim}{\longrightarrow}\Gamma(\tilde{S},\cO_{\tilde{S}}).
\end{equation}

\subsection{ }
\label{intro-2}
In the present paper, following  \cite{KR},
we define asymptotic vertex algebras, in particular, asymptotic algebras of chiral differential operators (ACDO).  Note that if a scheme $M$ is Poisson, then the jet scheme
$J_\infty M$ is vertex Poisson, i.e., $\cO_{J_\infty M}$ is a vertex Poisson algebra; this is proved in \cite{A1} but the fact that  brackets on jet spaces are of a '`vertex nature''
has been known to specialists in infinite dimensional dynamical systems for decades, see the introduction to \cite{BD} and references therein. Typically, an ACDO is a quantization
of  this or that $J_\infty M$. In particular,  given an algebra of chiral differential operators (CDO) $\cD^{ch}_M$ over $M$, see \cite{BD,MSV,GMS}, there is a simple construction  
of an ACDO $\cD^{ch}_{M,\sqrt{\hbar}}$, which is a sheaf over $J_\infty T^*M$.

Thus, if $X$ is a flag manifold, then there is a unique $\cD^{ch}_X$ \cite{GMSII}, and we construct an ACDO $ \cD^{ch}_{X,\sqrt{\hbar}}$, a sheaf over
$J_\infty T^*X$.  From $ \cD^{ch}_{X,\sqrt{\hbar}}$, via a chiral Hamiltonian reduction, we obtain $ \cD^{ch}_{\tilde{S},\sqrt{\hbar}}$, a certain ACDO over $J_\infty\tilde{S}$.

Similarly, given the well-known  centerless vertex algebra $V^{-h^\vee}_0(\fg)$, we construct an asymptotic vertex algebra $V^{-h^\vee}_0(\fg)_{\sqrt{\hbar}}$, from which, again via
a chiral Hamiltonian reduction, we obtain an asymptotic version of the centerless affine W-algebra, $\cW^{-h^\vee}_{0,\sqrt{\hbar}}$.  

Both $ \cD^{ch}_{\tilde{S},\sqrt{\hbar}}$ and $\cW^{-h^\vee}_{0,\sqrt{\hbar}}$ carry a $\CC^*$-action; in particular, $(\cW^{-h^\vee}_{0,\sqrt{\hbar}})^{\CC^*}$ is the usual
affine W-algebra, $\cW^{-h^\vee}_{0}$, as defined in \cite{FF2} for a principal nilpotent element and in \cite{KRW} in general. Our main result is the following two isomorphisms
\begin{equation}
\label{main-res-intro}
\Gamma(J_\infty \tilde{S},\cD^{ch}_{\tilde{S},\sqrt{\hbar}})\stackrel{\sim}{\longrightarrow}\cW^{-h^\vee}_{0,\sqrt{\hbar}},\;
\Gamma(J_\infty \tilde{S},\cD^{ch}_{\tilde{S},\sqrt{\hbar}})^{\CC^*}\stackrel{\sim}{\longrightarrow}\cW^{-h^\vee}_{0}.
\end{equation}
Notice that there is an obvious jet scheme version of (\ref{konst-thm-first-app-for-slodow-intro})

\begin{equation}
\label{konst-thm-first-app-for-slodow-jet-intro}
\CC[J_\infty S]\longrightarrow\Gamma(J_\infty\tilde{S},\cO_{J_\infty\tilde{S}}),
\end{equation}
but it is by no means an isomorphism, which makes (\ref{main-res-intro}) less obvious and indicates that the quantum objects are better suited for regularization.

\subsection{ }
\label{intro-3}

A few more words about the contents of this paper. Sect.~\ref{alg-as-chi-diff-oper} is devoted to, so to say, basics of asymptotic vertex algebras and ACDO. In particular,
we show that an ACDO may quantize, apart from the standard $J_\infty T^*M$, their twisted versions. We call them, perhaps
clumsily, twisted jet-cotangent bundles, and their isomorphism classes are classified by
$H^1(M,\Omega^2_M\rightarrow\Omega^{3,cl}_M)$.  This construction is modeled on the construction of twisted cotangent bundles in \cite{BB2}, but note that the latter
are classified by the cohomology group of a lower degree: $H^1(M,\Omega^1_M\rightarrow\Omega^{2,cl}_M)$.

In sect.~\ref{hamiltonian-reduction-subsection} we deal with the chiral Hamiltonian reduction. While the ideas are the old ones of Feigin-Kostant/Sternberg-Feigin/Frenkel, we
find that when working over an arbitrary smooth variety (at one point we need to relax even the smoothness condition) the use of asymptotic objects becomes essential.
Note that the constructed via the chiral Hamiltonian reduction $\cD^{ch}_{\tilde{S},\sqrt{\hbar}}$ are the first examples of an ACDO that do not come directly from some
CDO (as they are not related to a cotangent bundle). 

This approach gives rise to a novel definition of an affine W-algebra, see sect.~\ref{Definitions-of-vertex-Poisson-W-algebras-finite-and-affine}.  It is similar to the original
Feigin-Frenkel definition 
(which only works for a nilpotent element
admitting a good even grading)
and more intuitive than, but equivalent to, the general Kac-Roan-Wakimoto definition \cite{KRW}.

The main result can be found in sect.~\ref{we-form-main-res}. In sect.~\ref{Example-of-sl3} we explicitly work out the simplest nontrivial example that arises
when $\fg=sl_3$.
\subsection{Acknowledgments.} 
This work started when TA was visiting USC.
Part of the work on this was accomplished when FM visited IHES and MPIM, Bonn. 
TK is also grateful to MPIM for hospitality.

\section{algebras of asymptotic chiral differential operators}
\label{alg-as-chi-diff-oper}

By a vector space we will mean a vector superspace and by $\tilde{a}$ we will denote the parity of $a$.

\subsection{Vertex Poisson algebras and twisted jet-cotangent bundles}
\label{vert-poiss-alg}

\subsubsection{ }
\label{def-vert-poiss}
A¡¡{\em  vertex Poisson algebra}
 is a collection $(V,\mathbf{1}, \partial, _{(n)},n\geq -1)$, where 
 
 \begin{tabular}{l}
 $V$ is a  $\CC$-vector space, \\
$_{(n)}$ is a product $V\otimes V\rightarrow V$, $a\otimes b\mapsto a_{(n)}b$, \\
$\mathbf{1}\in V$ is a distinguished, {\em vacuum}, vector, \\
$\partial\in\text{End}_\CC V$,
\end{tabular}

that satisfies the following conditions:

(1) the collection $(V,_{(-1)},\vacA,\partial)$ is a commutative associative unital algebra with derivation $\partial$;

(2) the collection $(V,\partial,_{(n)},n\geq 0)$ is a vertex Lie algebra, i.e., the following holds

\begin{equation}
\label{loc-nilp-ax-poiss}
\text{ for all }a,b\in V,\; a_{(n)}b=0\text{ if }n\gg 0,
\end{equation}

\begin{equation}
\label{skew-symmetry}
a_{(n)}b=(-1)^{n+1+\tilde{a}\tilde{b}}\sum_{j=0}^{+\infty}\frac{(-1)^j}{j!}\partial^j(b_{(n+j)}a),
\end{equation}

 \begin{equation}
   \label{Borcherds-identity-poiss}
    a\ops{m} b\ops{k}c -
   (-1)^{\tilde{a}\tilde{b}}b\ops{k} a\ops{m}c =   \sum\limits_{j\geq 0} {m \choose j} (a\ops{j} b )\ops{m+k-j}\text{ for all } a,b,c\in V, n\geq 0,
 \end{equation}
 
\begin{equation}
\label{partial-derivation}
[\partial,a_{(n)}]=(\partial a)_{(n)}=-na_{(n-1)}\text{ for all }n\geq 0;
\end{equation}

(3) each $a_{(n)}$ is a derivation of the product $_{(-1)}$.

\bigskip

In what follows we will often unburden the notation by writing $ab$ instead of $a_{(-1)}b$.

\subsubsection{ }
\label{loc-vert-poiss-alg}
Let $V$ be a vertex Poisson algebra and $S\subset V$ a multiplicative  subset.  The derivation
$\partial$ extends  to the localization $ V[S^{-1}]$ by the `` quotient rule'':
\[
\partial(\frac{a}{s})=\frac{\partial(a)s-a\partial(s)}{s^2}.
\]
The $n$-th multiplication
\[
_{(n)}:V\otimes V\rightarrow V
\]
naturally extends to a map of localizations
\[
_{(n)}: V[S^{-1}]\otimes V[S^{-1}] \rightarrow V[S^{-1}].
\]
Indeed, the  $n$-th multiplication being a derivation, we define
\[
_{(n)}: V\otimes V[S^{-1}]\rightarrow V[S^{-1}]
\]
to be (the quotient rule again)
\[
x_{(n)}\frac{a}{s}=\frac{(x_{(n)}a)s-a(x_{(n)}s)}{s^2}.
\]
Next, we use  skew-symmetricity  to define
\[
_{(n)}:  V[S^{-1}]\otimes V\rightarrow V[S^{-1}]
\]
by the formula
\[
\frac{b}{t}_{(n)}a=(-1)^{n+1}\sum_{j=0}^{\infty}\frac{1}{j!}\partial^j(a_{(n+j)}\frac{b}{s}).
\]
Finally, we set
\[
\frac{b}{t}_{(n)}\frac{a}{s}=\frac{(\frac{b}{t}_{(n)}a)s-a(\frac{b}{t}_{(n)}s)}{s^2}
\]
\begin{lem}
\label{lem-on-loc-v-p-alg}
If $V$ is graded vertex Poisson algebra and $S\subset V$ is a multiplicative
subset, then the formulas written above define a vertex Poisson algebra on $V[S^{-1}]$.
\end{lem}
\begin{sloppypar}
{\em Proof.}  The map $_{(n)}:V\otimes V\rightarrow V$, $a\otimes b\mapsto a_{(n)}b$ is a (bi)differential operator: it is a derivation
in $b$ by definition, and the skew-symmetricity $a_{(n)}b=(-1)^{n+1}\sum_{j\geq 0}(-1)^j1/j!\partial^j(b_{(n+j)}a)$ shows that it is
a potentially infinite order differential operator in $a$. However, as follows from the same formula, the restriction of  $_{(n)}$ to any finitely generated subalgebra
of $V$ is of finite order, hence allows a canonical and standard extension to the localization  via the following recurrent
relation: if $P$ is a finite order differential operator, then define 
\end{sloppypar}
\begin{equation}
\label{local-of-diff-oper}
P(\frac{a}{s})=\frac{P(a)}{s}-\frac{1}{s}[P,s](\frac{a}{s}),
\end{equation}
noting that the order of $[P,s]$ is by definition less than that of $P$.

The formulas written before the lemma are but an explicit computation of these extensions.

To conclude,  notice that the vertex Poisson algebra axioms are then statements that certain differential operators are equal to zero, the property preserved upon
localization. \hfill $\qed$

\begin{cor}
\label{sheaf-of-vert-poiss-alg-over-spec}
If $V$ is a  vertex Poisson algebra, then the structure sheaf of the affine scheme $\text{Spec}V$ is a sheaf of vertex
Poisson algebras.
\end{cor}

\subsubsection{ }
\label{jet-schemes}
One can say that the spectrum of a vertex Poisson algebra is a vertex Poisson scheme, i.e., a scheme whose structure sheaf is a sheaf of vertex Poisson
algebras. Here is the avenue for producing not necessarily affine vertex Poisson schemes.

Let $X$ be a scheme, $J_\infty X$ the corresponding jet scheme. It is a $D$-scheme \cite{BD} and, in particular, the structure sheaf
$\cO_{J_\infty X}$ carries a canonical derivation $\partial$.

\begin{lem}\cite{A1}
\label{from-poiss-to-vert-poiss}
If $X$ is a Poisson  scheme, then $\cO_{J_\infty X}$ carries a unique vertex Poisson algebra structure such that, for each
$f,g\in\cO_X\subset \cO_{J_\infty X}$, $f_{(0)}g=\{f,g\}$ and $f_{(n)}g=0$ if $n>0$.
\end{lem}

\begin{rem}
\label{coisson-interpr}
The fact that the notion of a vertex Poisson algebra structure is a jet version of that of a Poisson structure, implicit in a 
body of work on infinite dimensional
dynamical systems, has been most eloquently put forward in \cite{BD} (see also references therein).  In the language of coisson algebras, proposed in \cite{BD},
the above lemma becomes practically obvious. The coisson algebra attached to a vertex Poisson algebra $V$ is $\cV=V\otimes\cO_\CC$, a $D$-module
on $\CC$ with coisson bracket
\[
\{.,.\}_{coiss}:\cV\boxtimes\cV\rightarrow \Delta_*\cV,\text{ s.t. } \{a(z),b(w)\}_{coiss}=\sum_{n=0}^{\infty}\frac{1}{j!}(a_{(n)}b)(w)\delta^{(n)}(z-w)
\text{ if }a,b\in V\otimes 1.
\]
A Poisson bracket on a commutative associative algebra $A$ in this language is reflected in letting
\[
\{a(z),b(w)\}_{coiss}=\{a,b\}(w)\delta(z-w).
\]
This easily and unambiguously extends to the jet-algebra $J_\infty A$, the $D$-algebra freely generated by $A$,
by repeatedly using the rule
\[
\{\partial_z^ia(z),\partial_w^jb(w)\}_{coiss}=\partial_z^i\partial_w^j(\{a(z),b(w)\}_{coiss}),
\]
which, by the way, allows one  quickly to write down somewhat complicated explicit formulas for $_{(n)}$-multiplications 
recorded in \cite{A1}.
\end{rem} 

\begin{cor}
\label{jinfty-is a-funct}
The assignment $X\mapsto J_\infty X$ defines a functor from the category of Poisson schemes to the category of vertex Poisson schemes.
\end{cor}

\subsubsection{ }
\label{coisson-on-dual-to-lie}
One typical example  is provided by a Lie algebra $\fm$. The dual $\fm^*$ is a Poisson variety (Kostant-Kirillov bracket),
$J_\infty\fm^*$ is a vertex Poisson scheme with $\CC[J_\infty \fm^*]=S^\bullet(\fm[t^{-1}]t^{-1})$. The meaning of the latter identification is that
a $\CC$-point of   $J_\infty\fm^*$  $f(t)\in\fm^*[[t]]$, where the value of a coordinate function $\phi(t)\in \fm[t^{-1}]t^{-1}$ is the residue $\text{Res}_{t=0}\phi(t)f(t)dt$.

\subsubsection{ }
\label{twist-jet-cot-bundle}
Perhaps the most obvious example of a vertex Poisson scheme suggested by Lemma~\ref{from-poiss-to-vert-poiss}
is $J_\infty T^*X$, $T^*X$ being equipped with the standard Poisson structure.  The vertex Poisson structure affords a very simple local description. 
If $\{x_j,\partial_j,1\leq j\leq \text{dim}X\}$
is a ``coordinate system'' on a sufficiently small $U\subset X$,
 meaning that $x_j\in\cO_X(U)$, $\partial_j\in\cT_X(U)$, and $\partial_i(x_j)=\delta_{ij}$, 
then $\partial_{i(0)}x_j=\delta_{ij}$ and all the products $_{(n)}$
with $n>0$ of $x_j$ and $\partial_j$ are zero. This assignment uniquely extends to the  algebra with derivation $\partial$ freely
 generated by $\{x_j,\partial_j,1\leq j\leq \text{dim}X\}$ and then   to the entire $\cO_{J_\infty T^*X}(J_\infty T^*U)$.
 Corollary~\ref{sheaf-of-vert-poiss-alg-over-spec} gives then a sheaf of vertex Poisson algebras on $\cO_{J_\infty T^*U}$.

This construction can be twisted. Consider a cover $\{U_s\}$ of $X$ with a coordinate system
 $\{x_j^{(s)},\partial_j^{(s)},1\leq j\leq \text{dim}X\}$ on each $U_s$ and a \v Cech cocycle
 $\{\alpha_s,\beta_{st}\}$  with coefficients in the complex $\Omega^2_X\rightarrow\Omega^{3,cl}_X$; this means that
$\alpha_s\in\Omega^{3,cl}_X(U_s)$, $\beta_{st}\in\Omega^2_X(U_s\cap U_t)$ and the following cocycle conditions hold:
\[
d\beta_{st}=\alpha_t-\alpha_s\text{ on }U_s\cap U_t,\; \beta_{st}-\beta_{rt}+\beta_{st}=0\text{ on } U_r\cap U_s\cap U_t.
\]

Now change the previously constructed vertex Poisson algebra structure on $\cO_{J_\infty T^*U_s}$ as follows:  recall that there is a natural embedding
$\Omega^1_X\stackrel{\sim}{\rightarrow}\cO_X\cdot\partial\cO_X\subset \cO_{J_\infty T^* X}$ for any $X$, declare  that
$\partial^{(s)}_{i(0)}\partial^{(s)}_j=\iota_{\partial^{(s)}_j}\iota_{\partial^{(s)}_i} \alpha_s\in\Omega^1_{U_s}$ , instead of 0, and
keep the other relations. Extended as above, this gives a new vertex Poisson algebra structure on
$\cO_{J_\infty T^* U_s}$. The reason the axioms hold true is the closedness of $\alpha_s$, an observation of Bressler \cite{Bre}.

The local change of the  vertex Poisson structure requires the change of the gluing functions. Over
 the intersections $U_s\cap U_t$, we re-glue by demanding that
$\partial^{(s)}_i\mapsto \partial^{(s)}_i+\iota_{\partial^{(s)}_i}\beta_{st}$. An easy to check formula 
\begin{equation}
\label{effect-add-2-form}
(\xi+\iota_{\xi}\beta)_{(0)}(\eta+\iota_{\eta}\beta)=\xi_{(0)}\eta+\iota_{[\xi,\eta]}\beta+\iota_{\eta}\iota_{\xi}d\beta
\end{equation}
implies that in this way we obtain a vertex Poisson algebra isomorphism
 \[
 \rho_{st}:\;\cO_{J_\infty T^* U_s}|_{J_\infty T^*U_s\cap J_\infty T^* U_t}\rightarrow \cO_{J_\infty T^* U_t}|_{J_\infty T^*U_s\cap J_\infty T^* U_t}.
 \]
That $\{\alpha_s,\beta_{st}\}$ is a cocycle implies the associativity $\rho_{rt}=\rho_{st}\circ\rho_{rs}$.
 Thus a sheaf of vertex Poisson algebras over $X$ arises; denote
it temporarily by $\cO_{J_\infty T^*X}^{tw}$. 

Notice that our manipulations have done  nothing to the associative commutative algebra structure locally.
 Therefore $\cO_{J_\infty T^*X}^{tw}(J\infty T^*U_{s})$ is (noncanonically) isomorphic
to $\cO_{J_\infty T^*X}(J_\infty T^*U_{s})$ and contains $\cO_X(U_s)$. This gives a morphism of schemes $\text{Spec}\cO_{J_\infty T^*X}^{tw}(J_\infty T^*U_s)
\rightarrow U_s$.

The gluing functions have been changed, but not their restriction to $\cO_X$. What this means is that the
 vertex Poisson schemes $\{\text{Spec}\cO_{J_\infty T^*X}^{tw}(U_s)\rightarrow U_s\}$
glue into vertex Poisson  scheme over $X$ that is locally, but not globally, isomorphic to $J_\infty T^*X$ as an $X$-scheme. 
{\em We denote this scheme by} $J_\infty^{tw} T^*X$ {\em and call
it a twisted jet-cotangent bundle over $X$.}

If we want to emphasize that $J_\infty^{tw} T^*X$ is built using a particular cocycle $\{\alpha_s,\beta_{st}\}$, then instead of $J_\infty^{tw} T^*X$ 
we shall write $J_\infty^{\alpha,\beta} T^*X$ 

It is rather clear that although the construction of $J_\infty^{\alpha,\beta} T^*X$ depends on the choice of a cocycle $\{\alpha_s,\beta_{st}\}$, 
the isomorphism class of the result
depends only on the cohomology class $\overline{\{\alpha_s,\beta_{st}\}}\in H^{1}(X,\Omega^2_X\rightarrow\Omega^{3,cl}_X)$.

\begin{rem}
\label{polish-way-say}
Our discussion is modeled on \cite{BB2}, 2.1.5--2.1.8, where twisted cotangent bundles are introduced.
It is easy to see that all twisted jet-cotangent bundles form a category, in fact a groupoid. Each such bundle is an $\Omega^2_X\rightarrow\Omega^{3,cl}_X$-torsor, 
and the category of  twisted jet-cotangent bundles 
is equivalent to that of $\Omega^2_X \rightarrow \Omega^{3,cl}_X$-torsors. We leave the details out.
\end{rem}
\subsubsection{ }
\label{super-jet-coiss-lambda}
There is little doubt that much of what we have said and will say has superanalogues. Here is a simple example of what we shall actually need.
If $M$ is a smooth supermanifold, then $T^*M$ is a symplectic supermanifold, and $J_\infty T^*M$ is a vertex super-Poisson scheme. For example,
if $\fm$ is a purely even vector space, then by $\Pi\fm$ we shall denote the corresponding purely odd vector space; $T^*\Pi\fm$ is a symplectic 
supermanifold, $J_\infty T^*\Pi\fm$ is a vertex super-Poisson scheme so that $\CC[J_\infty T^*\Pi\fm]$ is a vertex Poisson (super)algebra. We have
an identification $\CC[J_\infty T^*\Pi\fm]=S^\bullet(\fm[t^{-1}]t^{-1}\oplus\fm^*[t^{-1}]t^{-1})$ analogous  to the one discussed in sect.~\ref{coisson-on-dual-to-lie}.

It is clear that $\CC[T^*\Pi\fm]$ is nothing but the exterior algebra $\Lambda(\fm\oplus\fm^*)$ (with an obvious Poisson bracket), and we
shall denote by $\Lambda^{vert}(\fm\oplus\fm^*)$ the vertex Poisson algebra $\CC[J_\infty T^*\Pi\fm]$. All of this is well known, of course.

\subsection{Asymptotic CDO}
\label{asympt-cdo}
Now we wish to quantize $J_\infty^{tw} T^*X$. One expects this to mean constructing a sheaf of vertex algebras over $J_\infty^{tw} T^*X$  whose quasiclassical limit 
is
the structure sheaf $\cO_{J_\infty^{tw} T^*X}$. In essence, such a quantization is known -- a CDO  $\cD^{ch}_X$, \cite{MSV,GMS}. 
However, a CDO is a 
rather ``noncommutative'' object, and although it localizes along $X$, it resists localization over the entire
  $J_\infty^{tw} T^*X$. A way out was suggested by Kashiwara and
Rouquier \cite{KR} in the context of usual algebras of differential operators, and we will use it as a template.

\subsubsection{ }
\label{def-of asympt-vert-alg}
Recall  that a {\em  vertex algebra}
 is a collection $(V,\mathbf{1}, _{(n)},n\in\ZZ)$, where 
 
 \begin{tabular}{l}
 $V$ is a  $\CC$-vector space, \\
$_{(n)}$ is a product $V\otimes V\rightarrow V$, $a\otimes b\mapsto a_{(n)}b$, \\
$\mathbf{1}\in V$ is a distinguished, {\em vacuum}, vector, 
\end{tabular}

that satisfies the following conditions:
\begin{equation}
\label{vacuum-axiom}
a_{(-1)}\vacA=a,\; a_{(n)}\vacA=0\text{ for all }a\in V,\; n\geq 0,
\end{equation}
\begin{equation}
\label{loc-nilp-ax}
\text{ for all }a,b\in V,\; a_{(n)}b=0\text{ if }n\gg 0,
\end{equation}

 \begin{align}
   \label{Borcherds-identity}
   &\sum\limits_{j\geq 0} {m \choose j} (a\ops{n+j} b )\ops{m+k-j}\\
   =& \sum\limits_{j\geq 0} (-1)^{j} {n \choose j}\{ a\ops{m+n-j} b\ops{k+j} - (-1)^{n+\tilde{a}\tilde{b}}
   b\ops{n+k-j} a\ops{m+j} \}.
   \nonumber
   \end{align}
\bigskip   

Note that this definition suppresses the following important part of structure: the map
\[
\partial: V\rightarrow V,\; a\mapsto a_{(-2)}\vacA
\]
is a derivation of all products and satisfies (\ref{partial-derivation}).

\bigskip

 Following \cite{Li1}, we will say that an {\em $\hbar$-adic vertex algebra}  is a collection $(V,\mathbf{1}, _{(n)},n\in\ZZ)$, where 
 
 \begin{tabular}{l}
 $V$ is a flat $\CC[[\hbar]]$-module complete in $\hbar$-adic topology, \\
 $_{(n)} :V\otimes V\rightarrow V$ is a continuous $\CC[[\hbar]]$-linear map,\\
\end{tabular}

 such that for each $N\geq 0$ the collection $(V/\hbar^N V,\mathbf{1}, _{(n)},n\in\ZZ)$ is a vertex algebra.
 
\bigskip
Finally, and here we loosely follow the terminology suggested in \cite{DK}, a collection $(V,\mathbf{1}, _{(n)},n\in\ZZ)$ is called
 an {\em asymptotic} vertex algebra if $V$ is a $\CC((\hbar))$-vector space and there is an $\hbar$-adic vertex algebra  $(V',\mathbf{1}', _{(n)'},n\in\ZZ)$
 along with an embedding of $\CC[[\hbar]]$-modules $V'\rightarrow V$ so that $\vacA'\mapsto\vacA$, $V=V'\otimes_{\CC[[\hbar]]}\CC((h))$, and
  each multiplication $_{(n)}$ is an extension of
   of the multiplication $_{(n)'}$  by $\CC((\hbar))$-linearity. Any such $V'\subset V$ is called a {\em lattice} of $V$.
 
 \bigskip
 Here is a rather trivial example: if $V$ is a vertex algebra, then $V[[\hbar]]$ and $V((\hbar))$ with operations extended by
 $\CC[[\hbar]]$-  and $\CC((\hbar))$-linearity (resp.)   are an $\hbar$-adic and asymptotic 
 vertex algebra respectively.
 Less trivial examples
 will appear in sects.~\ref{defn-and-hadization-of-g},~\ref{defn-and-hadization-of-Lambda},~\ref{hadication-of-cdo}, and~\ref{main-def}
\begin{rem}
\label{asymtp-vert-vs-vert}
  Notice that although the notion of a vertex algebra makes sense over any ring that contains $\QQ$, an $\hbar$-adic vertex algebra is
not a vertex algebra over $\CC[[\hbar]]$, just as an asymptotic vertex algebra is not  a vertex algebra over $\CC((\hbar))$,
because the local nilpotency condition  (\ref{loc-nilp-ax}) does not quite hold; it fails even in the trivial example just considered.
  Typically and more generally,  we will have an $\hbar$-adic vertex algebra
$V=W[[\hbar]]$ for some $\CC$-vector space $W$, and for each $a,b\in W$
\[
a_{(n)}b=\sum_{k=0}^{\infty}a_{(n,k)}b\hbar^k, \;a_{(n,k)}b \in W.
\]
The definition as stated implies not that $  a_{(n)}b=0$ if $n\gg 0$ but that, for each $k$,  $a_{(n,k)}b=0$ if $n\gg 0$.

The Borcherds identity (\ref{Borcherds-identity}) holds on the nose and, in fact, vertex algebra techniques tend easily to apply to these
more general situations by focusing on the coefficients of appropriate $\hbar$-expansions. 
\end{rem}

\subsubsection{ }
\label{quasiclass-lim-hadic-vert}
If $V$ is an $\hbar$-adic vertex algebra, then $\hbar V\subset V$ is an ideal and
$V/\hbar V$ is a vertex algebra over $\CC$. If, furthermore, $V/\hbar V$ is a commutative vertex algebra,
i.e., $a_{(n)}b \in \hbar V$ for all $a,b\in V$, $n\geq 0$, then the standard definition
\[
(a\text{ mod }(\hbar))_{(n)}(b\text{ mod }(\hbar))=(\frac{1}{\hbar} a_{(n)}b)
\text{ mod }(\hbar)\text{ if }n\geq
0
\]
makes $V/\hbar V$ into a vertex Poisson algebra.

As usual, we will call an $\hbar$-adic vertex algebra $V$ a {\em quantization} of a vertex Poisson algebra
$W$ if $V/\hbar V$ is commutative and isomorphic to $W$ as a vertex Poisson algebra. For the lack of a better
name, we will also call an asymptotic vertex algebra $V$ a {\em quantization} of a vertex Poisson algebra
$W$ if one of its lattices is. The same definitions
apply to sheaves of $\hbar$-adic vertex algebras.

\subsubsection{ }
\label{defn-and-hadization-of-g}
Let $\fg$ be a finite dimensional Lie algebra with invariant inner product $(.,.)$.  A well-known vertex algebra attached to it,
$V^{(.,.)}(\fg)$, can be defined as a vertex algebra strongly generated by $\fg$ with relations
\[
x_{(0)}y=[x,y],\; x_{(1)}y=(x,y)\vacA,\; x_{(n)}y=0\text{ if }n>1.
\]
As a notational matter, if one $(.,.)$ is understood fixed, then $V^{k}(\fg)$ will stand for $V^{k(.,.)}(\fg)$. If only zero $(.,.)$ exists, say, when
$\fg$ is nilpotent, then we shall write simply $V(\fg)$.

Here is a way to attach an $\hbar$-adic and asymptotic vertex algebra to $V^{(.,.)}(\fg)$.  Inside $V^{(.,.)}(\fg)[\hbar]$ consider a vertex subalgebra
generated by $\hbar\fg$. Equivalently, writing $\hat{x}$ for $\hbar x$, the above relations become 
\[
\hat{x}_{(0)}\hat{y}=\hbar\widehat{[x,y]},\; \hat{x}_{(1)}\hat{y}=\hbar^2(x,y)\vacA,
\]
and we denote, tentatively, by $\widehat{V}^{(.,.)}(\fg)$ the vertex algebra generated by these elements and relations over $\CC[\hbar]$. Finally,
let $V^{(.,.)}(\fg)_{\hbar,\geq}$ the $\hbar$-adic completion of $\widehat{V}^{(.,.)}(\fg)$. This is an $\hbar$-adic vertex algebra. Localizing by $\hbar^{-1}$
we obtain an asymptotic vertex algebra $V^{(.,.)}(\fg)_{\hbar}$ of which $V^{(.,.)}(\fg)_{\hbar,\geq}$ is a lattice.

It is immediate to deduce from the relations above that $V^{(.,.)}(\fg)_{\hbar,\geq}/\hbar V^{(.,.)}(\fg)_{\hbar,\geq}$ is canonically isomorphic to, hence a quantization of,
$S(\fm[t^{-1}]t^{-1})$, the vertex Poisson algebra introduced in sect.~\ref{coisson-on-dual-to-lie}.

Notice that there is an obvious multiplication preserving morphism
\[
V^{(.,.)}(\fg)\longrightarrow V^{(.,.)}(\fg)_{\hbar},\; x\mapsto \hbar^{-1}\hat{x}.
\]

Later we shall need an ``algebraic extension''  $V^{(.,.)}(\fg)_{\sqrt{\hbar}}=V^{(.,.)}(\fg)_{\hbar}\otimes_{\CC((\hbar))}\CC((\sqrt{\hbar}))$.

\subsubsection{ }
\label{defn-and-hadization-of-Lambda}
Similarly, given a vector space $\fm$, let us define a vertex Clifford vertex algebra, $Cl^{vert}(\fm\oplus\fm^*)$, to be a vertex algebra strongly generated by
$\fm\oplus\fm^*$ with relations
\[
x_{(0)}y=\langle x,y\rangle,\;  x_{(n)}y=0\text{ if }n>1, x,y\in\fm\oplus\fm^*,
\]
where $\langle.,.\rangle$ is the natural inner product on $\fm\oplus\fm^*$.

Along the lines of sect.~\ref{defn-and-hadization-of-g}, one obtains, by replacing $x\in\fm$ with $\hat{x}=\hbar x$, $y\in\fm^*$ with $\hat{y}=y$, and then
completing by positive powers of $\hbar$, an $\hbar$-adic vertex algebra $Cl^{vert}(\fm\oplus\fm^*)_{\hbar,\geq}$, and then the asymptotic vertex algebra
$Cl^{vert}(\fm\oplus\fm^*)_{\hbar}$.

It is immediate to see that $Cl^{vert}(\fm\oplus\fm^*)_{\hbar,\geq}$ is a quantization of $\Lambda^{vert}(\fm\oplus\fm^*)$ introduced in sect.~\ref{super-jet-coiss-lambda}.
Note also an obvious ``dividing out $\hbar$'' morphism
\[
Cl^{vert}(\fm\oplus\fm^*)\longrightarrow Cl^{vert}(\fm\oplus\fm^*)_{\hbar}.
\]

\subsubsection{ }
\label{an intro-form-to-cdo}
Recall that a vertex algebra $V$ is called $\ZZ_+$-graded if $V=V_0\oplus V_1\oplus V_2\oplus\cdots$
so that $(V_n)_{(i)}(V_m)\subset V_{n+m-i-1}$ for all $i$, where we set $V_n=0$ if $n<0$. The same definition
applies to sheaves of vertex algebras on a variety and one easily verifies  \cite{GMS}, for any such sheaf
$\cV=\cV_0\oplus \cV_1\oplus \cV_2\oplus\cdots$, that

\begin{tabular}{l}
$(\cV_0,_{(-1)},\vacA)$ is an unital associative commutative algebra;\\
$(\cV_0)_{(-1)}(\partial\cV_0)$  and $\cV_1/(\cV_0)_{(-1)}(\partial\cV_0)$
are  $\cV_0$-modules, the action being defined by $_{(-1)}$;\\
$\partial:\cV_0\rightarrow (\cV_0)_{(-1)}(\partial\cV_0)$ is a derivation;\\
$\cV_1/(\cV_0)_{(-1)}(\partial\cV_0)$ is a $\cV_0$-algebroid Lie,  $_{(0)}$ defining the Lie bracket on $\cV_1/(\cV_0)_{(-1)}(\partial\cV_0)$ and\\
 the action of $\cV_1/(\cV_0)_{(-1)}(\partial\cV_0)$ on $\cV_0$ by derivations.
\end{tabular}

\bigskip
Let $X$ be a smooth algebraic variety.
\begin{defn}
\label{def-cdo}
An algebra of chiral differential operators (CDO) is a sheaf of $\ZZ_+$-graded vertex algebras over
$X$ , $\cD^{ch}_X$, so that

(i) the triple $(\cD^{ch}_X)_0$ is isomorphic to $\cO_X$ as a unital associative commutative algebra;

(ii)  $((\cD^{ch}_X)_0)_{(-1)}((\partial(\cD^{ch}_X)_0))=\Omega^1_X$ as an  $\cO_X$-module with derivation;

(iii) the pair $(\cD^{ch}_X)_1/((\cD^{ch}_X)_0)_{(-1)}(\partial(\cD^{ch}_X)_0)$ is isomorphic to $\cT_X$ as
 an $\cO_X$-algebroid Lie;

(iv) $\cD^{ch}_X$ is strongly generated by $(\cD^{ch}_X)_0\oplus(\cD^{ch}_X)_1$.
\end{defn}

CDOs were introduced in \cite{GMS} (and independently and in a much more general setting of chiral algebras
 in \cite{BD}) and classified as
follows.
\begin{thm}\cite{GMS}
\label{thn-on-class-cdo}
(i) A CDO over $X$ exists iff $ch_2(\cT_X)=0$.

(ii) If $ch_2(\cT_X)=0$, then the isomorphism classes of CDOs over $X$ form a torsor over
$H^{1}(X,\Omega^2_X\rightarrow\Omega^{3,cl}_X)$.
\end{thm}

\subsubsection{ }
\label{cdo-as-a-practical-matter}
Let $ch_2(\cT_X)=0$ and we fix one  CDO over $X$, $\cD^{ch}_X$. As a practical matter
 (we refer the reader to \cite{GMS} for further details), this amounts
to having chosen the following data:

{\em local:}
a fine enough affine cover $\{U_s\}$ of $X$, an abelian $\cO_X$-basis
$\{\partial_1,\ldots,\partial_N\}$ of
$\cT_X(U_s)$, and a closed 3-form $\alpha_s\in\Omega^{3,cl}_X(U_s)$ for each $s$;

{\em gluing:} an $n$-tuple of 1-forms $\{\phi_1^{st},\ldots,\phi_n^{st}\}$ 
on each double intersection $U_s\cap U_t$.

\bigskip
The local data determine the vertex algebra $\cD^{ch}_X(U_s)$ as follows:

$\cD^{ch}_X(U_s)$ contains and is generated by the $\CC$-vector space
$\cO_X(U_s)\oplus(\oplus_j\CC\partial_j)$; furthermore, $\cO_X(U_s)\ni 1=\vacA\in 
\cD^{ch}_X(U_s)$;

the  products $_{(n)}$, $n\geq 0$, amongst generating elements are
\begin{equation}
\label{relations}
(\partial_j)_{(0)}(\partial_k)=\iota_{\partial_k}\iota_{\partial_j}\alpha,
(\partial_j)_{(0)}f=\partial_j(f),\;(\partial_j)_{(1)}(\partial_k)_{(1)}=f_{(n)}g=0,\; 
f,g\in\cO_X(U_\bullet).
\end{equation}
Here and elsewhere we often omit the superscripts  that indicate the relevant charts.

To describe $\cD^{ch}_X(U_s)$ as a vector space, consider the
canonical composite projection
\[
\pi: \;J_\infty T^*X\rightarrow T^*X\rightarrow X.
\]
There is a vector space isomorphism
\begin{equation}
\label{symb-vert-cdo-local}
\cO_{J_\infty T^*X}(\pi^{-1}U_s)\rightarrow \cD^{ch}_X(U_s)
\end{equation}
defined by
\begin{eqnarray}
&\partial^{p_1}(\partial_{i_1})\cdots\partial^{p_l}(\partial_{i_l})
\partial^{q_1}(f_1)\cdots\partial^{q_l}(f_l)\nonumber\\
&\mapsto
(\partial^{p_1}(\partial_{i_1}))_{(-1)}\cdots(\partial^{p_l}(\partial_{i_l}))_{(-1)}
(\partial^{q_1}(f_1))_{(-1)}\cdots(\partial^{q_l}(f_l))_{(-1)}\vacA,\nonumber
\end{eqnarray}
where $f_j\in\cO_X$, $\partial$ stands for the canonical derivation of both the left  and right hand
sides, and the products on the right, which are not associative, are nested on the right.

\bigskip
The gluing data determine the patching over the double intersection $U_s\cap U_t$  as follows:
\begin{equation}
\label{how-to-glue}
f\mapsto f,\; \partial_j\mapsto\partial_j+\phi_j^{st}.
\end{equation}

Taking all of this for granted, one can easily verify
\begin{lem}
\label{elem-prop-cdo}

(i) If $\{x_j,\partial_j\}$ is a coordinate system and $\alpha=0$, then the subalgebra of
$\cD^{ch}_X(U_\bullet)$ generated by $\{x_j,\partial_j\}$ is a $\beta\gamma$-system, i.e., the only
non-trivial products of generators marked by nonnegative integers   are $(\partial_i)_{(0)}x_j=\delta_{ij}$.

(ii) The pieces $\cO_X(U_\bullet)$, $\cO_X(U_\bullet)\partial\cO_X(U_\bullet)$ patch into
the subsheaves $\cO_X,\Omega^1_X\subset \cD^{ch}_X$.

(iii) The pieces $\oplus_j\cO_X(U_\bullet)\partial_j\oplus \cO_X(U_\bullet)\partial\cO_X(U_\bullet)$ patch into
the subsheaf $\cV\subset \cD^{ch}_X$ that fits into an exact sequence
\begin{equation}
\label{def-of-vert-algebroid}
0\rightarrow \Omega^1_X\rightarrow\cV\stackrel{\sigma}{\rightarrow}\cT_X\rightarrow 0.
\end{equation}
(iv) {\em Change of splitting}: if $\beta\in\Omega^{2}_X(U_\bullet)$, then
\[
(\partial_i+\iota_{\partial_i}\beta)_{(1)}(\iota_{\partial_j+\partial_j}\beta)=0,\;
(\partial_i+\iota_{\partial_i}\beta)_{(0)}(\partial_j+\iota_{\partial_j}\beta)=
\iota_{\partial_j}\iota_{\partial_i}\alpha+\iota_{\partial_j}\iota_{\partial_i}d\beta.
\]

\end{lem}

$\cV$ is known as a {\em vertex algebroid}; the subsheaf $\cO_X\oplus\cV$ generates $\cD^{ch}_X$,
 and technically $\cD^{ch}_X$ is defined in \cite{GMS}
 as the vertex enveloping
algebra of $\cV$. 

\subsubsection{ }

\label{hadication-of-cdo}
We wish to find a reasonable $\hbar$-adic version of $\cD^{ch}_X$.
Consider $\cD^{ch}_X[[\hbar]]$; it is a sheaf of $\hbar$-adic vertex algebras, 
see sect.~\ref{def-of asympt-vert-alg}, and it contains $\cV[[\hbar]]$ so that, cf. 
(\ref{def-of-vert-algebroid}),
\[
0\rightarrow \Omega^1_X[[\hbar]]\rightarrow\cV[[\hbar]]
\stackrel{\sigma}{\rightarrow}\cT_X[[\hbar]]\rightarrow 0,
\]
is exact. 

\bigskip
{\em As an intermediate object, introduce  $\widehat{\cD}^{ch}_X$, the $\hbar$-adic vertex subalgebra of $\cD^{ch}_X[[\hbar]]$
generated by $\cO_X$ and  $\sigma^{-1}(\cT_X[[\hbar]]\hbar)$.}  
 
\bigskip
It is clear that $\widehat{\cD}^{ch}_X/\hbar\widehat{\cD}^{ch}_X$ is a commutative vertex algebra, hence
a vertex Poisson algebra, sect.~\ref{quasiclass-lim-hadic-vert}. We will now identify it. Consider the
canonical composite projection
\[
\pi: \;J_\infty T^*X\rightarrow T^*X\rightarrow X.
\]
\begin{lem}
\label{quasiclass-lim-hadic-cdo}
The vertex Poisson algebras $\pi_*\cO_{J_\infty T^*X}$ and $\widehat{\cD}^{ch}_X/\hbar\widehat{\cD}^{ch}_X$
are canonically isomorphic.
\end{lem}
{\em Proof.} By definition and as reviewed in sect.~\ref{cdo-as-a-practical-matter} in greater detail, locally $\widehat{\cD}^{ch}_X$
is generated by $\widehat{\partial}_j\stackrel{\text{def}}{=}\hbar\partial_j$ and $\cO_X(U_\bullet)$. 
The nonvanishing
 products $_{(n)}$, $n\geq 0$, among generators are
\[
(\widehat{\partial}_j)_{(0)}(\widehat{\partial}_k)=\hbar^2\iota_{\partial_k}\iota_{\partial_j}\alpha,\;
(\widehat{\partial}_j)_{(0)}f=\hbar\partial_j(f).
\]
The prescription of sect.~\ref{quasiclass-lim-hadic-vert} immediately recovers the ``canonical commutation
relations'' of $\cO_{J_\infty T^*X}$, cf. sect.~\ref{twist-jet-cot-bundle}. As to the patching, 
(\ref{how-to-glue}) gives
$
f\mapsto f,\; \widehat{\partial}_j\mapsto\widehat{\partial}_j+\hbar\phi_j,
$
which modulo $\hbar$ becomes
$
f\mapsto f,\; \widehat{\partial}_j\mapsto\widehat{\partial}_j,
$ as it should. \hfill$\qed$

\bigskip
Lemma~\ref{quasiclass-lim-hadic-cdo} defines, by adjunction, a sheaf morphism
\begin{equation}
\label{ defines, by adjunction, a sheaf morphism}
\pi^{-1}(\widehat{\cD}^{ch}_X/\hbar\widehat{\cD}^{ch}_X)\rightarrow \cO_{J_\infty T^*X}.
\end{equation}
This is not an isomorphism -- and so $\pi^{-1}(\widehat{\cD}^{ch}_X/\hbar\widehat{\cD}^{ch}_X)$ is not a quantization of 
 $ \cO_{J_\infty T^*X}$, cf. sect.~\ref{quasiclass-lim-hadic-vert} -- because $\pi^{-1}(\widehat{\cD}^{ch}_X)$  is
insensitive to localization along the fibers, although
$\pi^{-1}(\widehat{\cD}^{ch}_X/\hbar\widehat{\cD}^{ch}_X)(\pi^{-1}U)$ and $\cO_{J_\infty T^*X}(\pi^{-1}U)$
are isomorphic for each  $U\subset X$. We will now construct the desired localization, to be denoted 
$\Delta(\pi^{-1}\widehat{\cD}^{ch}_X)$.

\subsubsection{ }
\label{sect-on-quant-thm}

Localization of vertex Poisson algebras that we performed in sect.~\ref{loc-vert-poiss-alg} is a solution of an obvious
universal problem. By contrast, vertex algebras do not afford even a  definition of localization, not immediately at least. Our approach is then opportunistic by
 necessity:
we will define the desired localization of $\pi^{-1}(\widehat{\cD}^{ch}_X)$ to be a quantization 
(sect.~\ref{quasiclass-lim-hadic-vert}) of its quasiclassical limit, $\cO_{J_\infty T^*X}$.

We place ourselves in the situation of sect.~\ref{cdo-as-a-practical-matter}, where we have an affine $U\subset X$ with a coordinate system
$\{\partial_j, x_j\}$. Isomorphism (\ref{symb-vert-cdo-local}) lifts to a vector space isomorphism $\pi^{-1}\widehat{\cD}^{ch}_X(\pi^{-1}U)
\rightarrow\cO_{J_\infty T^*X}(\pi^{-1}U)[[\hbar]]$. There is no canonical such isomorphism, but in the presence of a coordinate system we can make a choice: 
as in the proof of Lemma~\ref{quasiclass-lim-hadic-cdo}, let $\widehat{\partial}_i=\hbar\partial_i$,  send

\begin{eqnarray}
 &(\partial^{p_1}(\widehat{\partial}_{i_1}))_{(-1)}\cdots(\partial^{p_l}(\widehat{\partial}_{i_l}))_{(-1)}
(\partial^{q_1}(f_1))_{(-1)}\cdots(\partial^{q_l}(f_l))_{(-1)}\vacA\nonumber\\
&\mapsto\partial^{p_1}(\partial_{i_1})\cdots\partial^{p_l}(\partial_{i_l})
\partial^{q_1}(f_1)\cdots\partial^{q_l}(f_l)\nonumber
\end{eqnarray}
and then extend using $\CC[[\hbar]]$-linearity and continuity.

Having made this choice, we identify 
\[
_{(n)}:\; \widehat{\cD}^{ch}_X(U)\otimes \widehat{\cD}^{ch}_X(U)\rightarrow \widehat{\cD}^{ch}_X(U)
\]
with
\[
_{(n)}:\; \cO_{J_\infty T^*X}(\pi^{-1}U)[[\hbar]]\otimes \cO_{J_\infty T^*X}(\pi^{-1}U)[[\hbar]]\rightarrow \cO_{J_\infty T^*X}(\pi^{-1}U)[[\hbar]].
\]
The latter is determined by its restriction
\[
_{(n)}:\; \cO_{J_\infty T^*X}(\pi^{-1}U)\otimes \cO_{J_\infty T^*X}(\pi^{-1}U)\rightarrow \cO_{J_\infty T^*X}(\pi^{-1}U)[[\hbar]].
\]
The following is the key observation:
\begin{lem}
\label{chir-oper-diff-oper}
For each $n\in\ZZ$, $k\in\ZZ_+$ there is a finite order (bi)differential operator $P_{nk}(.,.)\in D_{\pi^{-1}U\times \pi^{-1}U}$  such that the map
\[
_{(n)}:\; \cO_{J_\infty T^*X}(\pi^{-1}U)\otimes \cO_{J_\infty T^*X}(\pi^{-1}U)\rightarrow \cO_{J_\infty T^*X}(\pi^{-1}U)[[\hbar]].
\]
is
\[
a_{(n)}b=\sum_{k\geq 0}\hbar^k P_{nk}(a,b).
\]
\end{lem}
We postpone the proof of this lemma until sect.~\ref{proof-of-lemma-on-difgf-op} and quickly
 conclude the construction of $\Delta(\pi^{-1}\widehat{\cD}^{ch}_X)$.  For each $f\in\cO_{J_\infty T^*X}(\pi^{-1}U)$ we define
 \begin{equation}
 \label{def-of-mult-on-localized-sheaf} 
_{(n)}:\; \cO_{J_\infty T^*X}(\pi^{-1}U_f)\otimes \cO_{J_\infty T^*X}(\pi^{-1}U_f)\rightarrow \cO_{J_\infty T^*X}(\pi^{-1}U_f)[[\hbar]].
\end{equation}
by extending the differential operators of Lemma~\ref {chir-oper-diff-oper} to the indicated localization via (\ref{local-of-diff-oper}).
The Borcherds identity (\ref{Borcherds-identity}), being an identity satisfied
by differential operators, is preserved upon localization. Hence the assignment
\[
\pi^{-1}U_f\mapsto \cO_{J_\infty T^*X}(\pi^{-1}U_f)
\]
defines a sheaf of $\hbar$-adic vertex algebras over $\pi^{-1} U$ that localizes $\pi^{-1}\widehat{\cD}^{ch}_X(\pi^{-1}U)$. Denote it temporarily
by $\Delta(\pi^{-1}\widehat{\cD}^{ch}_X)_{\pi^{-1}U}$.

Furthermore, upon the performed identification 
 the patching of various $\widehat{\cD}^{ch}_X(U)$ into a sheaf over $X$ is also defined by means of differential operators, see 
sect.~\ref{cdo-as-a-practical-matter}, which means  that the same gluing functions patch various  $\Delta(\pi^{-1}\widehat{\cD}^{ch}_X)_{\pi^{-1}U}$  into a sheaf 
of $\hbar$-adiv vertex algebras over 
$J_\infty T^*X$.  This is the looked for $\Delta(\pi^{-1}\widehat{\cD}^{ch}_X)$.  By  construction, 
 $\Delta(\pi^{-1}\widehat{\cD}^{ch}_X)/\hbar\Delta(\pi^{-1}\widehat{\cD}^{ch}_X)$
is isomorphic to $\cO_{J_\infty T^*X}$ as a vertex Poisson algebra. This proves

\begin{thm}
\label{result-on-localiz-of-cdo}
The above constructed sheaf of $\hbar$-adic vertex algebras, $\Delta(\pi^{-1}\widehat{\cD}^{ch}_X)$, is a
 quantization of $\cO_{J_\infty T^*X}$.
\end{thm}

\subsubsection{ Proof of Lemma~\ref{chir-oper-diff-oper}.}  
\label{proof-of-lemma-on-difgf-op}
A convenient way to encode all products of two elements $A,B$
 of a
vertex algebra is provided by the OPE:
\[
A(z)B(w)\sim \sum_{n\in\ZZ}\frac{A_{(n)}B(w)}{(z-w)^{n+1}}.
\]
Now suppose in the situation of sect.~\ref{cdo-as-a-practical-matter} all 3-forms $\alpha=0$.
Then the vertex subalegbra of $\cD^{ch}_X(U)$ generated by the coordinates $x_j,\partial_j$ is 
described by Lemma~\ref{elem-prop-cdo}. This is what
V.Kac calls ``free field theory,'' \cite{K} p.88, and in this case the OPE is easily computed via Wick's theorem,
\cite{K} p.87, as the sum over ``contractions.'' We symbolically record it thus:
\[
A(z)B(w)\sim\sum_{\text{contr}}\frac{\text{contr}(A(z)B(w))}{(z-w)^{\phi(\text{contr})}}.
\]
A formal application of the Taylor theorem then gives
\[
\text{contr}(A(z)B(w))\sim\sum_{l=0}^{\infty}\frac{(z-w)^l}{l!}\partial_z^l\text{contr}(A(z)|_{z=w}B(w)),
\]
hence the OPE. The point of this quick reminder is that the contracted field $\text{contr}(A(z)B(w))$ is a
certain bidifferential operator applied to the pair of arguments $A$ and $B$ -- an elementary fact following
from the definition; we leave it as an exercise for the conscientious reader. Likewise, $\partial_z$ is a 
differential operator, in fact, a vector field equal to the infinite sum of expressions such as
$(\partial^{m+1}x)\partial/\partial(\partial^{m}x)$.

Furthermore, if we collect  the summands in Wick's theorem with a fixed number of contractions, then for the
$\hbar$-adic version, $\widehat{\cD}^{ch}_X(U)$,  we obtain
\[
A(z)B(w)=\sum_{k=0}^\infty\hbar^k\sum_{\text{$k$-fold contr}}\frac{\text{contr}(A(z)B(w))}{(z-w)^{\phi(\text{contr})}}.
\]
where contracted field is given by a differential operator of finite order.

The entire algebra $\cD^{ch}_X(U)$ is not quite generated by the indicated fields and is not a ``free theory'',
nevertheless it is rather clear that the same formula applies if one chooses the ordering as we did in
 (\ref{symb-vert-cdo-local}).
 
 Finally, to get rid of the restrictive condition that $\alpha=0$, pick a point $p\in U$ and restrict
 to the formal neighborhood of a point $p$, i.e., embed $\widehat{\cD}^{ch}_X(U)$ into $\widehat{\cD}^{ch}_X(U_{\text{form},p})$. The latter
 algebra is like the former except that  all power series in $x$, not only Taylor extensions of 
 regular on $U$
 functions, are allowed. All the  formulas recorded in sect.~\ref{cdo-as-a-practical-matter} make sense 
 upon this completion, hence
 $\widehat{\cD}^{ch}_X(U_{\text{form},p})$ is an $\hbar$-adic vertex algebra and 
 $\widehat{\cD}^{ch}_X(U)\hookrightarrow\widehat{\cD}^{ch}_X(U_{\text{form},p})$ is an $\hbar$-adic vertex algebra embedding.
 The advantage gained upon completing is that  in the bigger
 algebra
 $\alpha$ is exact. Picking $\beta$ s.t. $d\beta=\alpha$ and performing a change of variables
  $\partial_j\mapsto \partial-\iota_{\partial}\beta$, we get rid of $\alpha$, see Lemma~\ref{elem-prop-cdo} (iv).
   Therefore, the completed algebra is a (completed)
  free field theory, hence products in it are given by finite order differential operators, hence  their restrictions
  to $\widehat{\cD}^{ch}_X(U)$ are also.\hfill $\qed$
  
  \subsubsection{ }
  The above can be twisted so as to quantize the twisted jet-cotangent
  bundle $J_\infty^{tw} T^*X$, which we introduced in sect.\ref{twist-jet-cot-bundle}.
  Once again we place ourselves in the situation of sect.~\ref{cdo-as-a-practical-matter}, where we have a CDO
  $\cD^{ch}_X$ defined by means of a covering $\{U_s\}$ with a fixed coordinate system $\{x_i,\partial_i\}$ 
   and
  a closed 3-form $\bar{\alpha}_s$ on each $U_s$. Now given a \v Cech 2-cocycle $\{\alpha_{s},\beta_{st}\}$
  with coefficients in $\Omega^2_X\rightarrow\Omega^{3,cl}_X$, we construct a new CDO, 
  to be denoted $\cD^{ch}_X\boxplus(\alpha,\beta)$; this is a key to Theorem~\ref{thn-on-class-cdo}(ii).
  
  For this, tear $\cD^{ch}_X$ apart and consider a bunch of CDOs, $\{\cD^{ch}_{U_s}\}$,the restrictions of 
  $\cD^{ch}_X$ to each chart.
   Next, redefine the $_{(0)}$-multiplication on $\cD^{ch}_{U_s}$ by letting $(\partial_i)_{(0)}(\partial_j)=
  \iota_{\partial_j}\iota_{\partial_i}\bar{\alpha}_s+\iota_{\partial_j}\iota_{\partial_i}\alpha_s$  and 
  keeping all the other relations.
  This defines a new CDO over $U_s$, to be denoted by $\cD^{ch}_{U_s}\boxplus\alpha_s$. In order to glue 
  the pieces
  back together the gluing functions must be adjusted, and the required adjustment is provided by another
   component of our cocycle: over the intersection
  $U_{s}\cap U_t$ we compose the old transition functions with the transformation
  \[
  \partial_i\mapsto\partial_i+\iota_{\partial_i}\beta_{st}.
  \]
  The cocycle condition implies that these are vertex algebra isomorphisms (Lemma~\ref{elem-prop-cdo} (iv) is 
  instrumental here) and, furthermore,
  these isomorphisms define a new CDO, $\cD^{ch}_X\boxplus(\alpha,\beta)$.
  
  Now extend scalars to obtain $\cD^{ch}_X[\hbar,\hbar^{-1}]$ and consider 
  $\cD^{ch}_X[\hbar,\hbar^{-1}]\boxplus\hbar^{-1}(\alpha,\beta)$.
  It is generated by a vertex algebroid $\cV$ that fits into an exact sequence, cf. Lemma~\ref{elem-prop-cdo},
  \[
  0\rightarrow\Omega^1_X[\hbar,\hbar^{-1}]\rightarrow\cV\stackrel{\sigma}{\rightarrow}\cT_X[\hbar,\hbar^{-1}]\rightarrow 0.
  \]
  Consider the subsheaf of $\cD^{ch}_X[\hbar,\hbar^{-1}]\boxplus\hbar^{-1}(\alpha,\beta)$ generated, 
  as a vertex algebra, by
  $\cO_X$ and the preimage $\sigma^{-1}(\cT_X[\hbar]\hbar)$; this is a sheaf of vertex algebras over 
  $\CC[\hbar]$. 
  Denote by $\widehat{\cD}^{ch,\alpha,\beta}_X$ its completion by positive powers of $\hbar$. This is a sheaf 
  of $\hbar$-adic vertex algebras,
  a twist of $\widehat{\cD}^{ch}_X$ by the cocycle $\{\alpha_{st},\beta_s\}$.
  
  Now consider the twisted jet-cotangent bundle $J_\infty^{\alpha,\beta}T^*X$, sect.\ref{twist-jet-cot-bundle},
   and the composite projection
  \[
  \pi:J_\infty^{\alpha,\beta}T^*X\rightarrow T^*X\rightarrow X.
  \]
  It is immediate to obtain a vertex Poisson algebra isomorphism
  \begin{equation}
  \label{why-quant-twisted}
  \widehat{\cD}^{ch,\alpha,\beta}_X/\hbar\widehat{\cD}^{ch,\alpha,\beta}_X\stackrel{\sim}{\longrightarrow}\pi_*\cO_{J_\infty^{\alpha,\beta}T^*X}.
  \end{equation}
  This is analogous to o Lemma~\ref{quasiclass-lim-hadic-cdo} and the proof is also: because of the choices
   made -- and the coefficient $\hbar^{-1}$ in front of the cocycle -- we obtain the ``bracket''
  \[
  (\widehat{\partial}_i)_{(0)}(\widehat{\partial}_j)=
  \hbar\iota_{\partial_j}\iota_{\partial_i}\alpha_{s}+\hbar^2(\ldots),
  \]
  the transition functions
  \[
  \widehat{\partial}_i\mapsto\widehat{\partial}_i+\iota_{\partial_i}\beta_{st}+\hbar(\ldots),
  \]
  and it remains to take the coefficient of $\hbar$ in the former and the coefficient of $1=\hbar^0$ in the
  latter to get $\pi_*\cO_{J_\infty^{\alpha,\beta}T^*X}$.
  
  The discussion in sects.~\ref{sect-on-quant-thm} and \ref{proof-of-lemma-on-difgf-op} carries over to the
   present situation word for word. We
  record, somewhat informally,  the result as follows:
 \begin{thm}
\label{result-on-localiz-of-cdo-twisted}
 Multiplications and transition functions on $\pi^{-1}\widehat{\cD}^{ch,\alpha,\beta}_X$ are
 bidifferential operators and define by the formulas of sect.~\ref{sect-on-quant-thm} a 
 quantization of $\cO_{J_\infty^{\alpha,\beta} T^*X}$, to be denoted by 
 $\Delta(\pi^{-1}\widehat{\cD}^{ch,\alpha,\beta}_X)$.
\end{thm}

  \subsubsection{ }
  \label{main-def}
  We   will unburden the notation by introducing
  \[
  \cD^{ch}_{X,\hbar,\geq}=\Delta(\pi^{-1}\widehat{\cD}^{ch,\alpha,\beta}_X),\;
  \cD^{ch}_{X,\hbar}=\Delta(\pi^{-1}\widehat{\cD}^{ch,\alpha,\beta}_X)\otimes_{\CC[[\hbar]]}\CC((\hbar)),
  \]
  where we have suppressed the dependence on some of the data, such as the cocycle $\{\alpha_{s},\beta_{st}\}$.
  The former is a sheaf $\hbar$-adic vertex algebras over $J_\infty^{tw} T^*X$, the latter is a sheaf
  of asymptotic vertex algebras, sect.~\ref{def-of asympt-vert-alg} over $J_\infty^{tw} T^*X$. We would
  like to create a reasonable framework for these objects by proposing the following working definition:
  \begin{defn}
  \label{defn-of-hafic-asympt-cdo}
  Let   $Y$ be a vertex Poisson scheme admitting an open cover  $\{U_\alpha\}$, each $U_\alpha$ being
  isomorphic to $J_\infty V_\alpha$ for some  smooth symplectic affine variety $V_\alpha$ with a symplectomorphism
$V_\alpha\rightarrow T^* \CC^n$. 
  
 (i)  A sheaf $\cA$ of $\hbar$-adic  vertex algebras over $Y$ is called
  an $\hbar$-adic CDO  if there is a vertex Posson algebra isomorphism
  $\cA/\hbar\cA\stackrel{\sim}{\longrightarrow}\cO_Y$.
  
  (ii)A sheaf $\cA$ of asymptotic  vertex algebras over $Y$ is called
  an asymptotic CDO if it contains an $\hbar$-adic CDO $\cA_\geq$ 
  so that $\cA$ isomorphic to $\cA_\geq\otimes_{\CC[[\hbar]]}\CC((\hbar))$.
  $\cA_\geq$ is called a lattice of $\cA$.
  \end{defn}
  Thus $\cD^{ch}_{X,\hbar}$ is an asymptotic CDO (ACDO)  with lattice $\cD^{ch}_{X,\hbar,\geq}$, an $\hbar$-adic CDO.
  
  \subsubsection{ }
  \label{subs-on-filtr-acdo}
  Let $Y$ be a vertex Poisson scheme and
   $\cA$ be an asymptotic vertex algebra with lattice $\cA_\geq$ that is a quantization of $\cO_Y$.  Define a
decreasing filtration
\begin{equation}
\label{filtr-on-acdo}
\cA\supset\cdots \supset F_n\cA\supset F_{n+1}\cA\supset\cdots,\text{ where }F_n\cA=\{a\in\cA: \hbar^{-n}a\in\cA_\geq\}.
\end{equation}
One has 
\begin{equation}
\label{a-coll-ex-seq}
\cup_{n\ZZ}F_n\cA=\cA,\; \cap_{n\in\ZZ}F_n\cA=\{0\},\; F_n\cA/F_{n+1}\cA\stackrel{\sim}{\longrightarrow}\hbar^n\cO_{J_\infty Y}.
\end{equation}
Furthermore, the prescription of sect.~\ref{quasiclass-lim-hadic-vert} defines a vertex Poisson algebra structure on the graded
object, $\text{Gr}_F\cA$ and one immediately notices an isomorphism
\[
\text{Gr}_F\cA\stackrel{\sim}{\longrightarrow}\cO_{Y}\otimes\CC[\hbar,\hbar^{-1}],
\]    
where  $\CC[\hbar,\hbar^{-1}]$ is equipped with the obvious multiplication , $\partial$ and $_{(n)}$, $n\geq 0$ being 0.
\begin{defn}
\label{def-of-symbol}
Given an $a\in\cA$, let  $n$ be the greatest number such that $a\in F_n\cA$. Define the symbol
of $a$ to be
\[
\sigma(a)=a\text{ mod }F_{n+1}\cA\in \hbar^n\cO_{J_\infty Y}.
\]
\end{defn}
\subsubsection{ }
\label{emd-usu-cdo-into-asympt}
To construct an ACDO  $\cD^{ch}_{X,\hbar}$ we needed to be given a CDO $\cD^{ch}_X$. The two sheaves are closely related.
To formulate the result, denote by $G_n\cD^{ch}_X\subset \cD^{ch}_X$ the subsheaf consisting of those sections whose
local presentations as elements of $\cO_{J_\infty X}$, see isomorphism (\ref{symb-vert-cdo-local}), requires at most $n$ vector fields $\partial_i$.
It is clear that this defines an increasing filtration and the graded object is canonically isomorphic with $\cO_{J_\infty X}$.

Now recall that $\cD^{ch}_{X,\hbar}$ is a localization of (the pull-back of) $\widehat{\cD}^{ch}_X$, and the latter is generated 
by a part of $\cD^{ch}_X[\hbar]$; locally, the latter is generated by $\{\partial_i\}$, the former by $\{\widehat{\partial}_i=\hbar\partial_i\}$, cf. sect.\ref{hadication-of-cdo}. Hence
a tautological undoing morphism
\[
\cD^{ch}_X\rightarrow\widehat{\cD}^{ch}_X
\]
locally defined by $\partial_i\mapsto\hbar^{-1}\widehat{\partial}_i$.

\begin{lem}
\label{lem-emb-usu-cdo-asympt}
There is a morphism 
\[
\iota:\pi^{-1}\cD^{ch}_X\rightarrow \cD^{ch}_{X,\hbar},
\]
locally defined by $\partial_i\mapsto\hbar^{-1}\widehat{\partial}_i$. This morphism preserves all multiplications, $\iota(a)_{(n)}\iota(b)=\iota(a_{(n)}b)$,
and filtrations, $\iota(\pi^{-1}G_n\cD^{ch}_X)\subset F_{-n}\cD^{ch}_{X,\hbar}$.
\end{lem}

{\em Proof:} obvious.

\begin{cor}
\label{from-emd-in-cdo-to-emb-in-acdo}
A vertex algebra morphism $V\rightarrow\Gamma(X,\cD^{ch}_X)$ lifts to a morphism $V\rightarrow\Gamma(J_\infty T^*X,\cD^{ch}_{X,\hbar})$.
\end{cor}

 \subsection{Hamiltonian Reduction}
 
 \label{hamiltonian-reduction-subsection}

\subsubsection{ }
\label{ham-red-basic-setup}
Let $X$ be a smooth algebraic Poisson variety which carries a Poisson action of  a connected algebraic group $M$. By this we mean  that  the associated
Lie algebra morphism $\fm\stackrel{\text{def}}{=}\text{Lie}M\rightarrow \Gamma(X,\cT_X)$ lifts to a Lie algebra morphism
\begin{equation}
\label{class-moment-map}
\mu:\fm\rightarrow \Gamma(X,\cO_X)
\end{equation}
known as a dual {\em moment map}. Dualizing we obtain the moment map, a Poisson scheme morphism,
 \begin{equation}
 \label{class-mom-map-undual}
 \mu^*: X\rightarrow\fm^*.
 \end{equation}

For each character $\chi\in(\fm^*)^M$, the scheme-theoretic preimage $(\mu^*)^{-1}(\chi)\subset X$ carries an action
of  $M$.  We will assume once and for all that  $(\mu^*)^{-1}(\chi)$ is an $M$-torsor, by which we mean that at least locally the action affords a cross-section. Then we
obtain a  quotient  $(\mu^*)^{-1}(\chi)/M$ and a projection
\begin{equation}
\label{def-proj-muchi-mod-m}
p: (\mu^*)^{-1}(\chi)\longrightarrow (\mu^*)^{-1}(\chi)/M.
\end{equation}

If $\chi$ is a regular value for the restriction of $\mu^*$ to each symplectic leaf of $X$ and  $(\mu^*)^{-1}(\chi)/M$ is smooth, then
  $(\mu^*)^{-1}(\chi)/M$ is canonically
a smooth algebraic  Poisson variety, known as the Hamiltonian reduction of $X$, cf. \cite{Va}, Theorem~7.31.

\subsubsection{ }
\label{ham-red-basic-setup-coho-vers}

Kostant and Sternberg \cite{KS}, inspired by Feigin \cite{Feig}, have suggested the following cohomological
interpretation of the Hamiltonian reduction.  Consider $\Lambda(\fm\oplus\fm^*)$, the Poisson algebra introduced in sect.~\ref{super-jet-coiss-lambda},
and then  $\cO_X\otimes\Lambda(\fm\oplus\fm^*)$.
It is a (sheaf) Poisson algebra, and it contains a special and well-known globally defined element:
\begin{equation}
\label{def-of-reduct-diff}
d_{\chi}=\sum_i( \mu(m_i)-\chi(m_i))\otimes \phi_i^*-\frac{1}{2}\sum_{i,j,k}1\otimes c_{ij}^k\phi_k\phi^*_i\phi^*_j,
\end{equation}
where $\{m_i\}$ is a basis of $\fm$ regarded as our Lie algebra (hence purely even),  $\{\phi_i\}$ and $\{\phi_i^*\}$ are its copy and its dual inside 
 $\Lambda(\fm\oplus\fm^*)$, hence odd, and $\{c_{ij}^k\}$ are the structure constants relative to the chosen basis.  The operator 
 $\{d_\chi,.\}\in\text{End}(\cO_X\otimes\Lambda(\fm\oplus\fm^*))$ is a degree 1 differential and a derivation. Hence for each open $U\subset X$ the cohomology
 Poisson algebra $H^{\bullet}_{d_\chi}(\cO_X(U)\otimes\Lambda(\fm\oplus\fm^*))$ arises.  Denote by $\cH^{\bullet}_{d_\chi}(\cO_X\otimes\Lambda(\fm\oplus\fm^*) )$
 the sheaf associated to the presheaf $U\mapsto H^{\bullet}_{d_\chi}(\cO_X(U)\otimes\Lambda(\fm\oplus\fm^*))$.  We would like to compare two Poisson algebras,
 $p_*\cH^{0}_{d_\chi}(\cO_X\otimes\Lambda(\fm\oplus\fm^*) ) $, $p$ having been defined in (\ref{def-proj-muchi-mod-m}), and $\cO_{ (\mu^*)^{-1}(\chi)/M}$ that arises under
 the  extra regularity assumptions made in sect.~\ref{ham-red-basic-setup}.
 
 We shall say that the sequence  $\{\mu(m_i)-\chi(m_i)\}$ is {\em regular} if  for any $x\in (\mu^*)^{-1}(\chi)$ it defines a regular sequence in the local ring
 $\cO_{X,x}$.
 
 \begin{thm}
 \label{res-kost-sternb}
 (cf. \cite{KS}) If the sequence  $\{\mu(m_i)-\chi(m_i)\}$ is regular, then
 $\cH^{i}_{d_\chi}(\cO_X\otimes\Lambda(\fm\oplus\fm^*) )=0$  unless $i\geq 0$.  The sheaves $\cH^{i}_{d_\chi}(\cO_X\otimes\Lambda(\fm\oplus\fm^*) )$ 
 are supported on $(\mu^*)^{-1}(\chi) $ and  there is an algebra isomorphism
 \[
 p_*\cH^{0}_{d_\chi}(\cO_X\otimes\Lambda(\fm\oplus\fm^*) )\stackrel{\sim}{\longrightarrow}\cO_{ (\mu^*)^{-1}(\chi)/M}.
 \]
 Furthermore,  this is a Poisson algebra isomorphism provided  $\chi$ is a regular value for the restriction of $\mu^*$ to each symplectic leaf of $X$
and  $(\mu^*)^{-1}(\chi)/M$ is smooth.
 \end{thm} 
 
 \begin{rem}
 \label{kost-sternb-ext-class}
 Since $p_*\cH^{0}_{d_\chi}(\cO_X\otimes\Lambda(\fm\oplus\fm^*) )$ is automatically a Poisson algebra, this cohomological construction is a slight
 extension of the usual geometric one recalled in sect.~\ref{ham-red-basic-setup}.
 \end{rem}
 \bigskip
 Since the argument proving Theorem~\ref{def-of-reduct-diff} will reappear in one form or another in this note, we will now give a 
 
 {\em Proof:}  the question is local so we can assume that we are dealing with a $U\subset X$ so that the defining ideal of 
   $(\mu^*)^{-1}(\chi) \cap U$  is generated by the regular sequence $\{\mu(m_i)-\chi(m_i)\}$,
    and that the action of $M$ admits a cross-section $S$. The complex
   $(\cO_X(U)\otimes\Lambda(\fm\oplus\fm^*),\{d_\chi,.\})$ is given a bi-grading and a filtration by writing
\[
(\cO_X(U)\otimes\Lambda(\fm\oplus\fm^*))^n=\bigoplus_{p+q=n}\cO_X(U)\otimes\Lambda^p(\fm^*)\otimes\Lambda^{-q}(\fm),
\] 
\[
F^p(\cO_X(U)\otimes\Lambda(\fm\oplus\fm^*)^{p+q})=\bigoplus_{i\geq0}\cO_X(U)\otimes\Lambda^{p+i}(\fm^*)\otimes\Lambda^{-q+i}(\fm),
\]

Thus we obtain a spectral sequence $\{E^{pq}_r,d_r\}\Rightarrow H^{\bullet}_{d_\chi}(\cO_X(U)\otimes\Lambda(\fm\oplus\fm^*))$ such that
   
   $E_{0}^{pq}=K_q(\{\mu(m_i)-\chi(m_i)\}, \cO_X(U))\otimes\Lambda^p(\fm^*) $, where $K_q(\{\mu(m_i)-\chi(m_i)\}, \cO_X(U))$ is the Koszul complex of $\cO_X(U)$ w.r.t. the sequence  $\{\mu(m_i)-\chi(m_i)\}$;
   
   $E_1^{pq}=C^p(\fm,H_q^{Kosz}(\{\mu(m_i)-\chi(m_i)\}, \cO_X(U)))$ is the standard Chevalley complex of $\fm$ with coefficients in the Koszul cohomology;
   
   $E_2^{pq}=H^p(\fm,H_q^{Kosz}(\{\mu(m_i)-\chi(m_i)\}, \cO_X(U)))$ is then the cohomology of $\fm$ with the indicated coefficients.
   
It is now clear what happens:    If $(\mu^*)^{-1}(\chi)\cap U$ is empty, then the Koszul cohomology vanishes, and this proves the first assertion.
Otherwise, the sequence  $\{\mu(m_i)-\chi(m_i)\}$ being regular, the Koszul cohomology 
 equals  the algebra of functions on $(\mu^*)^{-1}(\chi)\cap U$ sitting in degree 0, which proves the second
 assertion. This also implies that the spectral sequence collapses in $\{E_2^{pq}\}$. 
 Furthermore, $E_2^{00}$ is the algebra of $\fm$-invariant functions, $\cO_{(\mu^*)^{-1}(\chi)}((\mu^*)^{-1}(\chi)\cap U)^{\fm}$, hence the  commutative associative
 algebra isomorphism of the theorem.  Checking  that this is also a Poisson algebra isomorphism 
  is routine and left to the reader.
  \hfill $\qed$
  
  \subsubsection{ }
  \label{jet-vers-ham-red-basic-setup-coho-vers} 
  Here is an obvious jet scheme version of sect.~ \ref{ham-red-basic-setup-coho-vers} obtained by repeatedly replacing
  various objects with their jet scheme versions, cf. Corollary~\ref{jinfty-is a-funct}.
   The moment map (\ref{class-mom-map-undual})
  engenders a map of vertex Poisson schemes
  \begin{equation}
  \label{mom-map-jet-vers-undual} 
  \mu^*_\infty:J_\infty X\rightarrow J_\infty\fm^*,
 \end{equation}
 hence a map of vertex Poisson algebras
 
  \begin{equation}
  \label{mom-map-jet-vers} 
  \mu_\infty:    S^\bullet(\fm\otimes\CC[t^{-1}]t^{-1})=\CC[J_\infty\fm^*]\rightarrow\Gamma(J_\infty X,\cO_{J_\infty X}).
 \end{equation} 
 
Consider the vertex Poisson algebra $\cO_{J_\infty X}\otimes\Lambda^{vert}(\fm\oplus\fm^*)$, where $\Lambda^{vert}(\fm\oplus\fm^*)$
is the vertex Poisson algebra we intoduced in sect.~\ref{super-jet-coiss-lambda}.  Note that since the projection
$J_\infty Y\rightarrow Y$ induces a canonical embedding $\cO_Y\hookrightarrow\cO_{J_\infty Y}$, the algebra 
$\cO_X\otimes\Lambda(\fm\oplus\fm^*)$
is canonically a subalgebra of $\cO_{J_\infty X}\otimes\Lambda^{vert}(\fm\oplus\fm^*)$, and
so the restriction of $\mu_\infty$ to $\fm$ coincides with $\mu$. Hence $\cO_{J_\infty X}\otimes\Lambda^{vert}(\fm\oplus\fm^*)$ acquires a globally
 defined element
$d_\chi$ defined by the familiar (\ref {def-of-reduct-diff}). This operator $(d_\chi)_{(0)}$ is a derivation ( a consequence of
(\ref{Borcherds-identity-poiss}) with $m=0$) and a differential, $(d_\chi)_{(0)}^2=0$.  Furthermore, $ \cO_{J_\infty X}\otimes\Lambda^{vert}(\fm\oplus\fm^*)$ is graded
by giving $\Pi\fm^*$ degree 1 and $\Pi\fm$ degree -1.
 Therefore, for each $U\subset J_\infty X$,
a cohomology vertex Poisson algebra, $H_{d_\chi}^{\infty/2+\bullet}(\cO_{J_\infty X}(U)\otimes
\Lambda^{vert}(\fm\otimes\fm^*))$, is defined.
Denote by $\cH_{d_\chi}^{\infty/2+\bullet}(\cO_{J_\infty X}\otimes
\Lambda^{vert}(\fm\otimes\fm^*))$ the sheaf associated with thus defined presheaf.

Projection (\ref{def-proj-muchi-mod-m}) determines
\[
p_\infty: J_\infty(\mu^*)^{-1}(\chi)\rightarrow J_\infty((\mu^*)^{-1}(\chi)/ M),
\]
which, under the assumptions of sect.~\ref{ham-red-basic-setup} is a  $J_\infty M$-torsor. 
Here is an analogue of
Theorem~\ref{res-kost-sternb}:
\begin{thm}
\label{anal-of-kost-sternb-puass} Let $X$ be a locally complete intersection algebraic variety,  $(\mu^*)^{-1}(\chi)$ be a reduced irreducible algebraic variety with rational singularities, and the sequence $\{\mu(m_i)-\chi(m_i)\}$ be regular. Then:

 (i) For each affine open $U\subset J_\infty X$, $H_{d_\chi}^{\infty/2+i}(\cO_{J_\infty X}(U)\otimes
\Lambda^{vert}(\fm\otimes\fm^*))=0$  unless $i\geq0$. 

(ii) For each affine open $U\subset J_\infty X$ s.t. $U\cap J_\infty(\mu^*)^{-1}(\chi)=\emptyset $, $H_{d_\chi}^{\infty/2+i}(\cO_{J_\infty X}(U)\otimes
\Lambda^{vert}(\fm\otimes\fm^*))=0$.
 
 (iii) $\cH_{d_\chi}^{\infty/2+i}(\cO_{J_\infty X}\otimes
\Lambda^{vert}(\fm\otimes\fm^*))=0$  unless $i\geq0$ and $\cH_{d_\chi}^{\infty/2+i}(\cO_{J_\infty X}\otimes
\Lambda^{vert}(\fm\otimes\fm^*))$ , $i\geq 0$,
 is supported on $J_\infty(\mu^*)^{-1}(\chi) $.

 (iv) If $U\subset J_\infty X$ be an affine open subset, then
 \[
 \cH_{d_\chi}^{\infty/2+0}(\cO_{J_\infty X}\otimes
\Lambda^{vert}(\fm\otimes\fm^*))(U)=H_{d_\chi}^{\infty/2+0}(\cO_{J_\infty X}(U)\otimes
\Lambda^{vert}(\fm\otimes\fm^*)).
\] 
 
 (v) There is a vertex Poisson algebra isomorphism
 \[
 p_{\infty*}\cH_{d_\chi}^{\infty/2+0}(\cO_{J_\infty X}\otimes
\Lambda^{vert}(\fm\otimes\fm^*))\stackrel{\sim}{\longrightarrow}\cO_{ J_\infty((\mu^*)^{-1}(\chi)/M)}.
 \]
\end{thm} 
\begin{rem}
\label{why-gen-lci} Note that  the familiar assumption that $X$ be smooth has been relaxed to include all locally complete
intersections. This generality will be needed when  dealing with $W$-algebras.
\end{rem}

\subsubsection{  Proof of Theorem~\ref{anal-of-kost-sternb-puass}.}  It is a straightforward variation of that of Theorem~\ref{res-kost-sternb}. We replace  $S$, $M$ with 
$J_\infty S$, $J_\infty M$ (resp.). The defining ideal of  $(\mu^*_\infty)^{-1}(\chi)\cap U$ is  generated by a collection of
functions that consists of  $\{\mu(m_i)-\chi(m_i)\}$ and all of its ``differential consequences'' $\{\partial^n(\mu(m_i)-\chi(m_i))\}$. Then Propositions 1.4, 1.5 and Theorem 3.3 of \cite{Mu} imply that the order $m$-jet scheme
 $J_m(\mu^*)^{-1}(\chi)$ is an irreducible and reduced complete intersection for each $m$. This implies that $\{\partial^n(\mu(m_i)-\chi(m_i))\}$
 is a regular sequence.

The complex $( \cO_{J_\infty X}\otimes\Lambda^{vert}(\fm\oplus\fm^*),( d_\chi)_{(0)})$ is bi-graded by assigning $\Pi\fm^*$ bi-degree $(1,0)$,
$\Pi\fm$ bi-degree $(0,-1)$.  Thus a spectral sequence  $\{E^{pq}_r(U),d_r\}$ arises with

 $E_{0}^{pq}(U)=K_q(\{\partial^n(\mu(m_i)-\chi(m_i))\}, \cO_{J_\infty X}(U))\otimes\Lambda^p(\fm^*[t^{-1}]t^{-1}) $,  the Koszul complex of $\cO_X(U)$ w.r.t. the sequence 
  $\{\partial^n(\mu(m_i)-\chi(m_i))\}$ tensored with    $\Lambda^p(\fm^*[t^{-1}]t^{-1})$;
   
   $E_1^{pq}(U)=C^p(\fm[t], H_q^{Kosz}(\{\partial^n(\mu(m_i)-\chi(m_i))\}, \cO_{J_\infty X}(U)))$ is the standard Chevalley complex of $\fm[t]$ with coefficients in the Koszul cohomology;
   
   $E_2^{pq}(U)=H^p(\fm[t],H_q^{Kosz}(\{\partial^n(\mu(m_i)-\chi(m_i))\}, \cO_{J_\infty X}(U)))$ is  the cohomology of $\fm[t]$ with the indicated coefficients.
   
   Note that in this infinite dimensional situation the bi-grading introduced is infinite and the convergence of the spectral sequence requires justification.
   Note that $\cO_{J_\infty X}\otimes\Lambda^{vert}(\fm\oplus\fm^*)$ is the structure sheaf  of the jet superscheme $J_\infty (X\times T^*\Pi\fm)$.
   The natural projection $J_\infty (X\times T^*\Pi\fm)\rightarrow X\times T^*\Pi\fm$ makes    $\cO_{J_\infty X}\otimes\Lambda^{vert}(\fm\oplus\fm^*)$ into a sheaf
   of $\cO_{X\times T^*\Pi\fm} $-modules. This sheaf is graded by ``the number of times the canonical derivation $\partial$ is applied'':
   \[
  \cO_{J_\infty X}\otimes\Lambda^{vert}(\fm\oplus\fm^*)=\bigoplus_{n=0}^{\infty} 
  \left(\cO_{J_\infty X}\otimes\Lambda^{vert}(\fm\oplus\fm^*)\right)_n.
  \] 
  It remains to notice that this direct sum decomposition is preserved by the differential, and each graded component has a bounded cohomological bi-grading.
  Hence $\{E^{pq}_r(U),d_r\}\Longrightarrow H_{d_\chi}^{\infty/2+p+q}(\cO_{J_\infty X}(U)\otimes
\Lambda^{vert}(\fm\otimes\fm^*))$.

\begin{sloppypar}
 The Koszul cohomology vanishes if
$(\mu^*_\infty)^{-1}(\chi)\cap U$ is empty and equals $\cO_{(\mu^*_\infty)^{-1}(\chi)}((\mu^*_\infty)^{-1}(\chi)\cap U)$ sitting in degree 0 otherwise. Therefore the
spectral sequence collapses and $E_2^{00}$ picks the algebra of $\fm[t]$-invariant functions. 
\end{sloppypar}
This proves items (i, ii,iii, v) of the theorem. Let us now do item (iv), a convenient computational tool.  Consider $\{\cE^{pq}_r,\tilde{d}_r\}$, a sheaf-theoretic
version of $\{E^{pq}_r(U),d_r\}$, where $\cE^{pq}_r$ is the sheaf associated to the presheaf $U\mapsto E^{pq}_r(U)$. For each $U$ there is an obvious
morphism of spectral sequences
\[
\{E^{pq}_r(U),d_r\}  \rightarrow \{\cE^{pq}_r(U),\tilde{d}_r\} .
\]
Furthermore, as we have seen, $E^{00}_1(U)=\cO_{(\mu^*_\infty)^{-1}(\chi)}((\mu^*_\infty)^{-1}(\chi)\cap U)$, therefore
 $\cE^{00}_1(U)=\cO_{(\mu^*_\infty)^{-1}(\chi)}((\mu^*_\infty)^{-1}(\chi)\cap U)$ too.  Hence an isomorphism
 \begin{equation}
 \label{proof-quas-comp-tool-1}
 E^{00}_2(U)\stackrel{\sim}{\longrightarrow}\cE^{00}_2(U).
 \end{equation}
 Now notice that there is a canonical morphism   
 \begin{equation}
 \label{proof-quas-comp-tool-2} 
 H_{d_\chi}^{\infty/2+0}(\cO_{J_\infty X}(U)\otimes
\Lambda^{vert}(\fm\otimes\fm^*))\longrightarrow \cH_{d_\chi}^{\infty/2+0}(\cO_{J_\infty X}\otimes
\Lambda^{vert}(\fm\otimes\fm^*))(U).
\end{equation}
It induces a morphism of the graded objects
\begin{equation}
\label{proof-quas-comp-tool-3}
 Gr H_{d_\chi}^{\infty/2+0}(\cO_{J_\infty X}(U)\otimes
\Lambda^{vert}(\fm\otimes\fm^*))\longrightarrow Gr\cH_{d_\chi}^{\infty/2+0}(\cO_{J_\infty X}\otimes
\Lambda^{vert}(\fm\otimes\fm^*))(U).
\end{equation}
As we have seen, $Gr H_{d_\chi}^{\infty/2+0}(\cO_{J_\infty X}(U)\otimes
\Lambda^{vert}(\fm\otimes\fm^*))=E^{00}_2(U)$ and $Gr \cH_{d_\chi}^{\infty/2+0}(\cO_{J_\infty X}\otimes
\Lambda^{vert}(\fm\otimes\fm^*))=\cE^{00}_2$.  Isomoprhism (\ref{proof-quas-comp-tool-1}) implies that
(\ref{proof-quas-comp-tool-3}) is an isomorphism, hence so is (\ref{proof-quas-comp-tool-2}).
\hfill $\qed$.

\subsubsection{ }
\label{quant-chital-hamilt-red}
We wish to quantize the discussion of sect.~\ref{jet-vers-ham-red-basic-setup-coho-vers}.    
The`` chiral Hamiltonian  reduction'' carried out in Lemma~\ref{anal-of-kost-sternb-puass} is a passage
from $\cO_{J_\infty X}$ to $\cO_{ J_\infty((\mu^*)^{-1}(\chi)/M)}$.  A quantization of the former, an asymptotic vertex algebra $\cA$, see 
sect.~\ref{quasiclass-lim-hadic-vert},
is assumed given. A quantization of the latter, which {\em a priori} does not have to exist, will be defined by quantizing the cohomology sheaf
$\cH_{d_\chi}^0(\cO_{J_\infty X}\otimes
\Lambda^{vert}(\fm\otimes\fm^*))$ or rather the complex $(\cO_{J_\infty X}\otimes
\Lambda^{vert}(\fm\otimes\fm^*), (d_\chi)_{(0)})$; this is known as the ``BRST approach to quantization.''  
A quantization of the vertex Poisson algebra $\Lambda^{vert}(\fm\otimes\fm^*)$ is  $Cl^{vert}(\fm\oplus\fm^*)_{\hbar}$,
see sect.~\ref{defn-and-hadization-of-Lambda}. Thus we  obtain $\cA\otimes Cl^{vert}(\fm\oplus\fm^*)$, where the tensoring is $\CC((\hbar))$-linear
and completed by positive powers of $\hbar$,
and what one really needs is a  chiral moment map, a quantization of (\ref{mom-map-jet-vers} ).

{\em One important assumption we need to make in order to proceed is that $M$ be unipotent.}

Let $V(\fm)$ be the vertex algebra attached to $\fm$, a well-known quantization of the vertex Poisson algebra  $S^\bullet(\fm\otimes\CC[t^{-1}]t^{-1})$, cf. 
sect.~\ref{coisson-on-dual-to-lie}; note
that $\fm$ being nilpotent no central extensions/charges are involved.  $V(\fm)$ carries a standard increasing  filtration, $\{G_nV(\fm)\}$,  so that
$Gr_GV(\fm)=S^\bullet(\fm\otimes\CC[t^{-1}]t^{-1})$. Given $v\in G_nV(\fm)\setminus G_{n-1}V(\fm)$, denote by $\sigma(v)\in S^n(\fm\otimes\CC[t^{-1}]t^{-1})$
its symbol.

$\cA$ also carries a  filtration,  albeit a decreasing one, $\{F_n\cA\}$, see sect.~\ref{subs-on-filtr-acdo}, hence so does $\Gamma(J_\infty X,\cA)$.

Given  map (\ref{mom-map-jet-vers} ), define a {\em chiral moment map} to be a vertex algebra morphism (cf. \cite{KR}, sect.2.4.1)
 \begin{equation}
  \label{mom-map-chiral} 
  \mu_{ch}:    V(\fm)\rightarrow\Gamma(J_\infty X,\cA)
 \end{equation} 
so that
\begin{equation}
\label{cond-on-symb-chir-mom}
\mu_{ch}(G_n V(\fm))\subset F_{-n}\Gamma(J_\infty X,\cA),\; \sigma(\mu_{ch}(v))=\mu_{\infty}(\sigma(v)),
\end{equation}
see Definition~\ref{def-of-symbol} in regard to symbol $\sigma(.)$ on the left. 

Consider the vertex algebra $\cA\otimes Cl^{vert}(\fm\oplus\fm^*)_{\hbar}$. It is graded by giving $\Pi\fm^*$ degree 1 and $\Pi\fm$ degree -1. Given a chiral 
moment map,
we obtain a globally defined element
\begin{equation}
\label{def-of-reduct-diff-chiral}
d_{\chi}^{ch}=\sum_i( \mu_{ch}(m_i)-\hbar^{-1}\chi(m_i))\otimes\hat{ \phi}_i^*-
\frac{1}{2}\sum_{i,j,k}\hbar^{-1}\otimes c_{ij}^k(\hat{\phi}_k)_{(-1)}((\hat{\phi}^*_i)_{(-1)}(\hat{\phi}^*_j)).
\end{equation}
The meaning of $\hat{\phi}$ and $\hat{\phi}^*$ is one explained in sect.~\ref{defn-and-hadization-of-Lambda}, and so this  element is an obvious quantization of
(\ref{def-of-reduct-diff}).

Note that since $m_i\in G_1V(\fm)$, $\mu_{ch}(m_i)-\hbar^{-1}\chi(m_i)\in F_{-1}\Gamma(J_\infty X,\cA)$.

The operator $(d_{\chi}^{ch})_{(0)}$ is a derivation and a degree 1 differential of $\cA\otimes Cl^{vert}(\fm\oplus\fm^*)$; that these properties quantize is also
a consequence of $\fm$ being nilpotent.
Denote by $\cH^{\infty/2+\bullet}_{d^{ch}_\chi}(\cA\otimes Cl^{vert}(\fm\oplus\fm^*)_{\hbar})$ the sheaf of asymptotic vertex algebras
 associated
with the presheaf  $U\mapsto H^{\infty/2+\bullet}_{d^{ch}_\chi}(\cA(U)\otimes Cl^{vert}(\fm\oplus\fm^*)_{\hbar})$.
\begin{thm}
\label{chiral-quant-hamilt-red}
Under the assumptions of  Theorem~\ref{anal-of-kost-sternb-puass}, 

(i) for each affine open $U\subset J_\infty X$,  $H^{\infty/2+i}_{d^{ch}_\chi}(\cA(U)\otimes Cl^{vert}(\fm\oplus\fm^*)_{\hbar})=0$ unless $i\geq 0$;

(ii) for each affine open $U\subset J_\infty X$ s.t.
$U\cap J_\infty(\mu^*)^{-1}(\chi)=\emptyset$, $H^{\infty/2+i}_{d^{ch}_\chi}(\cA(U)\otimes Cl^{vert}(\fm\oplus\fm^*)_{\hbar})=0$;

(iii) $\cH^{\infty/2+i}_{d^{ch}_\chi}(\cA\otimes Cl^{vert}(\fm\oplus\fm^*)_{\hbar})=0$ unless $i\geq 0$ and
$\cH^{\infty/2+i}_{d^{ch}_\chi}(\cA\otimes Cl^{vert}(\fm\oplus\fm^*)_{\hbar})$, $i\geq 0$, is supported on $J_\infty(\mu^*)^{-1}(\chi) $;

(iv) if $U\subset J_\infty X$ is an affine open subset, then
\[
\cH^{\infty/2+0}_{d^{ch}_\chi}(\cA\otimes Cl^{vert}(\fm\oplus\fm^*)_{\hbar})(U)=H^{\infty/2+0}_{d^{ch}_\chi}(\cA(U)\otimes Cl^{vert}(\fm\oplus\fm^*)_{\hbar});
\]

(v) $p_{\infty*}\cH^{\infty/2+0}_{d^{ch}_\chi}(\cA\otimes Cl^{vert}(\fm\oplus\fm^*)_{\hbar})$ is a quantization of $\cO_{J_\infty((\mu^*)^{-1}(\chi)/M)}$, and

$p_{\infty*}\cH^{\infty/2+i}_{d^{ch}_\chi}(\cA\otimes Cl^{vert}(\fm\oplus\fm^*)_{\hbar})=0$ if $i\neq0$.
\end{thm}
\subsubsection{ Proof of Theorem~\ref{chiral-quant-hamilt-red} (i,ii,iii).}
\label{proofofthmchquantham-i-ii-iii}

For each open $U\subset J_\infty X$,  $\cA(U)\otimes Cl^{vert}(\fm\oplus\fm^*)_{\hbar}$ is an asymptotic vertex algebra, hence it carries a filtration,
$F_\bullet(\cA(U)\otimes Cl^{vert}(\fm\oplus\fm^*)_{\hbar})$, by powers of $\hbar$, that we intoduced in sect.~\ref{subs-on-filtr-acdo}. This filtration is
obviously preserved by the differential, and so a spectral sequence, $\{E_r^{pq}(U), d_r\}$, arises. We have
\[
(E^{pq}_0(U), d_0)=\hbar^{p}( \cO_{J_\infty X}(U)\otimes\Lambda^{vert}(\fm\oplus\fm^*)^{p+q}, (d_{\chi})_{(0)}).
\]
This is simply because the whole construction is a quantization of the differential graded vertex Poisson algebra 
$( \cO_{J_\infty X}(U)\otimes\Lambda^{vert}(\fm\oplus\fm^*)^{p+q}, (d_{\chi})_{(0)})$. Therefore
\[
(E^{pq}_1(U), d_0)=\hbar^{p}H^{\infty/2+ p+q}_{d_\chi}( \cO_{J_\infty X}(U)\otimes\Lambda^{vert}(\fm\oplus\fm^*)).
\]

The spectral sequence does not really converge to $H^{\infty/2+\bullet}_{d^{ch}_\chi}(\cA(U)\otimes Cl^{vert}(\fm\oplus\fm^*)_{\hbar})$, but it does for some $U$
and in any case the $\hbar$-adic completeness of $\cA(U)\otimes Cl^{vert}(\fm\oplus\fm^*)_{\hbar}$ allows one  to 
derive properties of $H^{\infty/2+\bullet}_{d^{ch}_\chi}(\cA(U)\otimes Cl^{vert}(\fm\oplus\fm^*)_{\hbar})$ from Theorem~\ref{anal-of-kost-sternb-puass}.
For example, one has
\begin{equation}
\label {from-quasiclass-to-quant-1}
H^{\infty/2+i}_{d_\chi}( \cO_{J_\infty X}(U)\otimes\Lambda^{vert}(\fm\oplus\fm^*))=0\Longrightarrow 
H^{\infty/2+i}_{d^{ch}_\chi}(\cA(U)\otimes Cl^{vert}(\fm\oplus\fm^*)_{\hbar})=0.
\end{equation}

Indeed, assume given $a\in F_N(\cA(U)\otimes Cl^{vert}(\fm\oplus\fm^*)_{\hbar}^{i})$ that is annihilated by
$(d^{ch}_{\chi})_{(0)}$.  The assumption guarantees the existence  of  $b_0\in F_N(\cA(U)\otimes Cl^{vert}(\fm\oplus\fm^*)_{\hbar}^{i-1})$
s.t. $a-(d^{ch}_{\chi})_{(0)}b_0\in F_{N+1}(\cA(U)\otimes Cl^{vert}(\fm\oplus\fm^*)_{\hbar}^{i})$. An obvious induction will give then a sequence
$\{b_j\}$ s.t. $b_j\in F_{N+j}(\cA(U)\otimes Cl^{vert}(\fm\oplus\fm^*)_{\hbar}^{i-1})$ and 
\[
a-(d^{ch}_{\chi})_{(0)}(b_0+\cdots +b_m)\in F_{N+m+1}(\cA(U)\otimes Cl^{vert}(\fm\oplus\fm^*)_{\hbar}^{i}).
\]
Therefore
\[
a=(d^{ch}_{\chi})_{(0)}\sum_{j=0}^{\infty}b_j.
\]
Assertions (i) and (ii) now follow from the respective assertion of Theorem~\ref{anal-of-kost-sternb-puass}, and (iii) is a consequence of (i) and (ii). 
\subsection{ Proof of Theorem~\ref{chiral-quant-hamilt-red} (iv).}
\label{proofofthmchquantham-iv}
Along with $\{E^{pq}_r(U), d_r\}$ consider its sheaf version, $\{\cE^{pq}_r,\tilde{d}_r\}$, where $\cE^{pq}_r$ is the sheaf associated to the presheaf
$U\mapsto E^{pq}_r(U)$.  There is an obvious morphism of spectral sequences
\[
\{E^{pq}_r(U), d_r\}\rightarrow \{\cE^{pq}_r(U), \tilde{d}_r\}.
\]
Theorem~\ref{anal-of-kost-sternb-puass} (i) implies that
\[
E^{p,-p}_r(U)=\text{Ker}d_r\subset E^{p,-p}_1(U).
\]
Likewise
\[
\cE^{p,-p}_r(U)=\text{Ker}\tilde{d}_r\subset \cE^{p,-p}_1(U).
\]

Note that by Theorem~\ref{anal-of-kost-sternb-puass} (iv, v)
\[
 \oplus_p E^{p,-p}_1(U)= \cO_{J_\infty (\mu^*)^{-1}(\chi)/M}(U\cap J_\infty (\mu^*)^{-1}(\chi)/M)[\hbar,\hbar^{-1}], 
 \]
 hence
\[
\oplus_p \cE^{p,-p}_1(U)= \cO_{J_\infty (\mu^*)^{-1}(\chi)/M}(U\cap J_\infty (\mu^*)^{-1}(\chi)/M)[\hbar,\hbar^{-1}], 
\]
too. Therefore, the initial terms of the 2 spectral
sequences coincide; hence so do all of them. Thus we obtain:
\begin{equation}
\label{auxil-for-iv}
\cap_{r=1}^{\infty}\text{Ker}d_r=\cap_{j=1}^{\infty}\text{Ker}\tilde{d}_r.
\end{equation}
The relation of all of this to the problem at hand is as follows.  There is a natural morphism
\begin{equation}
\label{what-need-iv}
H^{\infty/2+0}_{d^{ch}_\chi}(\cA(U)\otimes Cl^{vert}(\fm\oplus\fm^*)_{\hbar})\longrightarrow\cH^{\infty/2+0}_{d^{ch}_\chi}(\cA\otimes Cl^{vert}(\fm\oplus\fm^*)_{\hbar})(U).
\end{equation}
It engenders a morphism of the graded objects
\begin{equation}
\label{what-need-iv-garded}
Gr H^{\infty/2+0}_{d^{ch}_\chi}(\cA(U)\otimes Cl^{vert}(\fm\oplus\fm^*)_{\hbar})\longrightarrow Gr\cH^{\infty/2+0}_{d^{ch}_\chi}(\cA\otimes Cl^{vert}(\fm\oplus\fm^*)_{\hbar})(U).
\end{equation}
Analogously to  (\ref{from-quasiclass-to-quant-1}) one verifies that
\[
H^{\infty/2+i-1}_{d_\chi}( \cO_{J_\infty X}(U)\otimes\Lambda^{vert}(\fm\oplus\fm^*))=0\Longrightarrow 
Gr H^{\infty/2+i}_{d^{ch}_\chi}(\cA(U)\otimes Cl^{vert}(\fm\oplus\fm^*)_{\hbar})=\cap_{r=1}^{\infty}\text{Ker}d_r.
\]
Likewise
\[
H^{\infty/2+i-1}_{d_\chi}( \cO_{J_\infty X}(U)\otimes\Lambda^{vert}(\fm\oplus\fm^*))=0\Longrightarrow 
Gr \cH^{\infty/2+i}_{d^{ch}_\chi}(\cA\otimes Cl^{vert} (\fm\oplus\fm^*)_{\hbar})=\cap_{r=1}^{\infty}\text{Ker}\tilde{d}_r.
\]
Thanks to the vanishing $H^{\infty/2-1}_{d_\chi}( \cO_{J_\infty X}(U)\otimes\Lambda^{vert}(\fm\oplus\fm^*))=0$, (\ref{auxil-for-iv}) implies that (\ref{what-need-iv-garded}) is an isomorphism, hence (\ref{what-need-iv}) is also, as desired.
\subsubsection{ Proof of Theorem~\ref{chiral-quant-hamilt-red} (v).}
\label{proofofthmchquantham-v}
Consider the projection $p_\infty:J_\infty(\mu^*)^{-1}(\chi)\rightarrow J_\infty(\mu^*)^{-1}(\chi)/J_\infty M$. We have assumed, see sect.~\ref{ham-red-basic-setup},
 that each point
of $J_\infty(\mu^*)^{-1}(\chi)/J_\infty M$ has a neighborhood, $V$, over which this projection  trivializes to become $V\times J_\infty M\rightarrow V$.
Pick $U\subset J_\infty X$ so that $U\cap J_\infty(\mu^*)^{-1}(\chi) \stackrel{\sim}{=}V\times J_\infty M$. For such $U$ our spectral sequence converges
and, in fact, collapses in $E^{pq}_1$.  Indeed, it follows from the proof of Theorem~\ref{anal-of-kost-sternb-puass} that
\[
E^{pq}_1(U)=\hbar^p H^{p+q}(\fm[t], \cO_{J_\infty X}(V\times J_\infty M)).
\]
By definition,
\[
H^{p+q}(\fm[t], \cO_{J_\infty X}(V\times J_\infty M))=\cO_{J_\infty(\mu^*)^{-1}(\chi)/J_\infty M}(V)\otimes H^{p+q}_{DR}(J_\infty M).
\]
Since $M$ is unipotent, the de Rham cohomology $H^{p+q}_{DR}(J_\infty M)=0$ unless $p+q=0$. Hence the  spectral sequence collapses. Its convergence is proved
by the inductive argument completely parallel to the  one used to prove (\ref{from-quasiclass-to-quant-1}). Hence an isomorphism (where Theorem~\ref{chiral-quant-hamilt-red} (iv) is used)
\[
Gr p_{\infty*}\cH^{\infty/2+0}_{d^{ch}_\chi}(\cA\otimes Cl^{vert}(\fm\oplus\fm^*)_{\hbar})\stackrel{\sim}{\longrightarrow }
\cO_{J_\infty(\mu^*)^{-1}(\chi)/J_\infty M}[\hbar,\hbar^{-1}],
\]
compatible by construction with the vertex algebra structure, and the vanishing of higher cohomology. This completes the proof. \hfill$\qed$

\section{$W$-algebras}
\label{sect-w-algebras}
Throughout this section, $\fg$ will denote a simple complex finite dimensional Lie algebra, $(.,.)$ a fixed once and for all
invariant inner product, hence an identification $\fg\stackrel{\sim}{\longrightarrow}\fg^*$. 
\subsection{Set-up}
\label{set-up-w-alg} 
\subsubsection{ }
\label{inside-g} 
Let a nonzero  $f\in\fg$ be nilpotent. Include it into an $sl_2$-triple, $\{f,h,e\}$, where $[e,f]=h$, $[h,e]=2e$, $[h,f]=-2f$.
Eigenvalues of $\text{ad}h$ define a grading $\fg=\oplus_i\fg_i$. The space $\fg_1$ carries a symplectic form 
\[
\omega: \fg_1\otimes\fg_1\rightarrow\CC,\; a\otimes b\mapsto (f,[a,b]).
\]
Choose a Lagrangian subspace $l\subset \fg_1$. Denote by $\fm=\fm_l$ the subspace $l\oplus(\bigoplus_{i\geq 2}\fg_i)\subset\fg$.
It is a nilpotent Lie subalgebra, and $\chi\stackrel{\text{def}}{=}(f,.)\in\fm^*$ is a character. Denote by $M$ the unipotent
group with Lie algebra $\fm$. 

One finds oneself in the situation of sect.\ref{ham-red-basic-setup-coho-vers}, and here we follow the remarkably lucid \cite{GG}:

$\fg^*$ is a Poisson manifold;

the (restriction of the) coadjoint action $M\times\fg^*\rightarrow\fg^*$ is Poisson;

the moment map is the obvious restriction map $\mu^*:\fg^*\rightarrow\fm^*$, the dual moment map is an embedding
$\mu:\fm\hookrightarrow\fg$;

$(\mu^*)^{-1}(\chi)=f+\fm^\perp$;

the $M$-space $f+\fm^\perp$ has the {\em Slodowy slice} $\SSS=f+\fg^e$ for a cross-section, i.e., the action of $M$ defines
an isomorphism of algebraic varieties  
\begin{equation}
\label{GG-Lemma}
M\times\SSS\rightarrow f+\fm^\perp.
\end{equation}

It follows that $\SSS$ is a Poisson smooth algebraic  variety. If $\cN$ is the nilpotent cone, then $S\stackrel{\text{def}}{=}\SSS\cap\cN$
is a Poisson subscheme of $\SSS$, in fact,  a reduced, irreducible, normal complete intersection \cite{Pre}.

Note an equivalent point of view: $\cN$ is a Poisson variety Poisson acted upon by $M$ of which $S$ is the Hamiltonian reduction w.r.t. character
$(f,.)$.
\subsubsection{ }
\label{set-up-w-alg-flag}
{\em For the rest of this note, $X$ will stand for the flag manifold of $\fg$.} The natural Lie algebra morphism $\fg\rightarrow \Gamma(X,\cT_X)$,
upon extending to an algebra morphism $S^\bullet\fg\rightarrow S^\bullet\Gamma(X,\cT_X)$, determines a morphism of algebraic varieties
\[
T^*X\mapsto \cN\subset\fg,
\]
known as the Springer resolution. Denote by $\tilde{S}$ the preimage of $S$ under this map.

Note that both $T^*X$ and $\cN$ are Poisson varieties carrying a Poisson action of $M$, and the Springer resolution obviously respects both
these structures. It follows that $\tilde{S}$ is  the Hamiltonian reduction of $T^*X$ w.r.t. $(f,.)$ and the natural map
\begin{equation}
\label{s-tilde-s}
\Phi:\tilde{S}\rightarrow S,
\end{equation}
is a symplectic a resolution , \cite{G}, Proposition 2.1.2.  Since $S$ is normal, we obtain a Poisson algebra isomorphism
\begin{equation}
\label{iso-of-funct-spaces}
\Phi^{\#}:\CC[S]\stackrel{\sim}{\longrightarrow}\Gamma(\tilde{S},\cO_{\tilde{S}}).
\end{equation}
\begin{lem}
\label{S-ratnl-sing}
 $S$ is an irreducible reduced rational singularities complete intersection scheme.
 \end{lem}
 {\em Proof:} 
 
 since (\ref{s-tilde-s}) is a symplectic resolution, it is a crepant resolution \cite{Fu}, Proposition~1.6;

existence of a crepant resolution implies $S$ has canonical singularities, by definition;

being a complete intersection, hence Gorenstein,  $S$ has canonical singularities if and only if it has rational singularities, cf. \cite{Mu}, discussion after Theorem~0.1;

that $S$ is a reduced, irreducible complete intersection is the above cited result of Premet. \hfill $\qed$

\subsection{W-algebras, finite and affine}
\label{Definitions-of-vertex-Poisson-W-algebras-finite-and-affine}
\subsubsection{ }
\label{w-ald-def-and-quant}
Considerations of sect.~\ref{inside-g} and \ref{set-up-w-alg-flag} give us a collection of Poisson algebras, $\CC[\fg^*]$, $\CC[\SSS]$,
$\CC[S]$, $\Gamma(\tilde{S}, \cO_{\tilde{S}})$. Hence (Lemma~\ref{from-poiss-to-vert-poiss}) a collection of vertex Poisson algebras
$\CC[J_\infty\fg^*]$, $\CC[J_\infty\SSS]$,
$\CC[J_\infty S]$, $\Gamma(J_\infty\tilde{S}, \cO_{J_\infty\tilde{S}})$. Let us discuss  quantizations of these 8 algebras.

Quantization of $\CC[\fg^*]$ is, of course, the universal enveloping algebra $U(\fg)$.  Quantization of $\CC[\SSS]$ was found  by Premet \cite{Pre}; 
this is a finite $W$-algebra attached to $\fg$ and $f\in\fg$ to be denoted $\cW^{fin}(\fg,f)$. (To be fair, the construction makes a choice of $l\subset\fg_1$, 
sect.~\ref{set-up-w-alg}, but the result is known to be independent of the choice \cite{GG}.) One possible definition \cite{GG} uses the Hecke algebra philosophy:
consider the generalized Gelfand-Graev module, $H_f=U(\fg)\otimes_{U(\fm)}\CC_\chi$, $\chi=(f,.)$, and then put $\cW^{fin}(\fg,f)=\text{End}_\fg^{opp} H_f$.

The latter definition suggests a simple way to quantize $\CC[S]$. We have an obvious algebra morphism $\cZ(\fg)\rightarrow \cW^{fin}(\fg,f)$ and for
any character $\nu\in\text{MaxSpec}\cZ(\fg)$ we let $\cW^{fin}_\nu(\fg,f)$ be the quotient of $\cW^{fin}(\fg,f)$ by the kernel of $\nu$. Each 
$\cW^{fin}_\nu(\fg,f)$ is a reasonable quantization of $\CC[S]$, especially $\cW^{fin}_0(\fg,f)$.

The problem of quantizing $\Gamma(\tilde{S}, \cO_{\tilde{S}})$ is what the machinery of asymptotic differential operators was invented for in \cite{KR}; we will
not be talking about it as we are focused    on its chiral analogues.

\subsubsection{ }
\label{chiral-quant-hadic}
To the jet scheme versions. We shall now propose a novel approach to affine $W$-algebras based on asymptotic 
vertex algebras to be compared later with the previously known one due to Kac, Roan, and Wakimoto \cite{KRW}.
As a quantization  of $\CC[J_\infty\fg^*]$, we choose $V^k(\fg)_{\sqrt{\hbar}}$, see sect.~\ref{defn-and-hadization-of-g}; note the appearance of
the central charge $k$. Now notice that according to sect.~\ref{jet-vers-ham-red-basic-setup-coho-vers}  $\CC[J_\infty\SSS]$ can be written as
$H_{d_\chi}^{\infty/2+0}(\CC[J_\infty\fg^*]\otimes
\Lambda^{vert}(\fm\otimes\fm^*))$, the 0th cohomology of the Hamiltonian reduction complex $ (\CC[J_\infty\fg^*]\otimes
\Lambda^{vert}(\fm\otimes\fm^*), d_\chi)$, $\chi=(f,.)$. Quantization of this complex is immediate and analogous to what we did in
sect.~\ref{quant-chital-hamilt-red}. Define the chiral Hamiltonian reduction complex to be
\begin{eqnarray}
\label{chir-ham-red-compl-usu-defn}
&(V^k(\fg)_{\sqrt{\hbar}}\hat{\otimes} Cl^{vert}(\fm\oplus\fm^*)_\hbar, (d^{ch}_{\chi})_{(0)}),\\
&d^{ch}_{\chi}= \sum_i\hbar^{-1}( m_i-(f,m_i))\otimes \phi_i^*-
\frac{1}{2}\sum_{i,j,k}1\otimes c_{ij}^k(\phi_k)_{(-1)}((\phi^*_i)_{(-1)}(\phi^*_j)),\nonumber
\end{eqnarray}
cf. (\ref{def-of-reduct-diff-chiral}).

 We then define  an affine asymptotic  $W$-algebra to be 
\[
\cW^k(\fg,f)_{\sqrt{\hbar}}=H_{d_{\chi}^{ch}}^{\infty/2+0}(V^k(\fg)_{\sqrt{\hbar}}\hat{\otimes}
Cl^{vert}(\fm\otimes\fm^*)_\hbar).
\]
Quantization of $\CC[J_\infty S]$ is dealt with analogously. If the central charge $k=-h^\vee$, then $V^k(\fg)$ acquires a big center $\fz(\fg)\subset V^k(\fg)$,
\cite{FF2}.  The vertex algebra $\fz(\fg)$ is $\ZZ_+$-graded, cf. sect.~\ref{an intro-form-to-cdo}:
$\fz(\fg)=\oplus_{n\geq 0}\fz(\fg)_n$ with $\fz(\fg)_0=\CC\vacA$. We  introduce $\fz(\fg)_+=\oplus_{n>0}\fz(\fg)_n$.  The subspace 
$(\fz(\fg)_+)_{(-1)}
V^{-h^\vee}(\fg)\subset V^{-h^\vee}(\fg)$ is a vertex ideal and we define a vertex algebra 
 $V^{-h^\vee}(\fg)_0=V^{-h^\vee}(\fg)/(\fz(\fg)_+)_{(-1)}V^{-h^\vee}(\fg)$. Likewise one obtains an asymptotic version $V_0^{-h^\vee}(\fg)_{\sqrt{\hbar}}$.
 Define
\[
\cW^{-h^\vee}_0(\fg,f)_{\sqrt{\hbar}}=H_{d_{\chi}^{ch}}^{\infty/2+0}(V^{-h^\vee}_0(\fg)_{\sqrt{\hbar}}\hat{\otimes}
Cl^{vert}(\fm\otimes\fm^*)_\hbar).
\]
\begin{lem}
\label{why-quantizs-of-jet-slice}

(i) $H_{d_{\chi}^{ch}}^{\infty/2+i}(V^k(\fg)_{\sqrt{\hbar}}\hat{\otimes}
Cl^{vert}(\fm\otimes\fm^*)_\hbar)=0$ if $i\neq 0$ and $\cW^k(\fg,f)_{\sqrt{\hbar}}$ is a quantization of $\CC[J_\infty\SSS]$.

(ii) $H_{d_{\chi}^{ch}}^{\infty/2+i}(V^{-h^\vee}_0(\fg)_{\sqrt{\hbar}}\hat{\otimes}
Cl^{vert}(\fm\otimes\fm^*)_\hbar)=0$ if $i\neq 0$ and $\cW^{-h^\vee}_0(\fg,f)_{\sqrt{\hbar}}$ is a quantization of $\CC[J_\infty S]$.
\end{lem}
{\em Proof.}  Since $\SSS$ is smooth and a cross-section for the action of $M$, see (\ref{GG-Lemma}), item (i) is an immediate consequence of Theorem~\ref{chiral-quant-hamilt-red}(v) with $X=\fg$. 
The same argument works in (ii) if we let  $X=\cN$, cf. the end of sect.~ \ref{inside-g}, thanks to Lemma~\ref{S-ratnl-sing}. \hfill $\qed$

\begin{rem}\label{rem:central-quotient}
As in   \cite{A3}, we see that
 the natural map
$\fz(\fg) \rightarrow \cW^{-h^{\vee}}(\fg,f)_{\sqrt{\hbar}}
$
is injective and that
$\cW_0(\fg,f)_{\sqrt{\hbar}}=\cW^{-h^{\vee}}(\fg,f)_{\sqrt{\hbar}}/
(\fz(\fg)_+)_{(-1)}
\cW^{-h^{\vee}}(\fg,f)_{\sqrt{\hbar}}$.
\end{rem}
\subsubsection{ }
\label{c-star-action-on-asympt-w and-affine-w}
In order to define vertex algebras $\cW^k(\fg,f)$ and $\cW^{-h^\vee}_0(\fg,f)$, i.e., to get rid of $\hbar$, one has to introduce a $\CC^*$-module structure
on the chiral Hamiltonian reduction complex.  For $t\in\CC^*$ consider the assignment
\begin{equation}
\label{beginn-of-jet-kazhd-gr-hadic}
\hat{x}_{(-1)}\vacA\mapsto t^2\hat{x}_{(-1)}\vacA,\; \psi\mapsto\psi,\; \hbar\mapsto t^2\hbar, x\in\fg,\psi\in(\fm\oplus\fm^*)\subset Cl^{vert}(\fm\oplus\fm^*).
\end{equation}
It is clear that this extends to an automorphism of $(V^k(\fg)_{\sqrt{\hbar}}\hat{\otimes} Cl^{vert}(\fm\oplus\fm^*)_\hbar)$ to be denoted
$\rho(t)$

Next, if $h$ is the element of the $sl_2$-triple introduced in sect.~\ref{set-up-w-alg}, then $ad_h$ is a diagonalizable derivation of $\fg$ that preserves $\fm\subset\fg$.
For this reason, $h$ extends to a diagonalizable derivation of  $(V^k(\fg)_{\sqrt{\hbar}}\hat{\otimes} Cl^{vert}(\fm\oplus\fm^*)_\hbar)$.  Furthermore, it is easy to see
that $t^{-h}\rho(t)$, unlike $\rho(t)$, respects the differential $d^{ch}_\chi$.

Define a group homomorphism
\[
\CC^*\longrightarrow \text{Aut}_{\CC}(V^k(\fg)_{\sqrt{\hbar}}\hat{\otimes} Cl^{vert}(\fm\oplus\fm^*)_\hbar, (d^{ch}_{\chi})_{(0)}),\; t\mapsto t^{-h}\rho(t).
\]
 All of this is but an extension of the celebrated Kazhdan filtration to the asymptotic vertex algebra situation, cf.\cite{GG}, but notice that what was a filtration has become
a  $\CC^*$-action thanks to $\hbar$. This
 is a beautiful idea of Kashiwara and Rouquier \cite{KR}.

Now introduce
\[
\cW^{k}(\fg,f)=\cW^{k}(\fg,f)_{\sqrt{\hbar}}^{\CC^*},\; \cW^{-h^\vee}_0(\fg,f)=\cW^{-h^\vee}_0(\fg,f)_{\sqrt{\hbar}}^{\CC^*}.
\]
{\em A priori} these are subalgebras of asymptotic vertex algebras and may contain infinite sums, but in fact one has:
\begin{lem}
\label{why-actually-vert-alg-1}
$\cW^{k}(\fg,f)=\cW^{k}(\fg,f)_{\sqrt{\hbar}}^{\CC^*}$ and $\cW^{-h^\vee}_0(\fg,f)=\cW^{-h^\vee}_0(\fg,f)_{\sqrt{\hbar}}^{\CC^*}$ are vertex algebras.
\end{lem}
{\em Proof.}  Lemma~\ref{why-quantizs-of-jet-slice} implies an isomorphism of $\CC^*$-modules
\[
\cW^k(\fg,f)_{\sqrt{\hbar}}\stackrel{\sim}{\longrightarrow}\CC[J_\infty\SSS]((\hbar)),
\]
where the action of $\CC^*$ on the r.h.s. is determined by the Kazhdan grading of $\CC[\SSS]$, see \cite{GG}, 4.1,4.2.  In fact, one defines the action of $\CC^*$ on
the entire $\CC[J_\infty\fg^*]$ by the same assignment $t\mapsto t^{-h}\rho(t)$ except that now it is stipulated that $\rho(t)x=t^2x$, $x\in\fg$; cf. (\ref{beginn-of-jet-kazhd-gr-hadic}). This action descends to $\CC[\SSS]$
and defines the Kazhdan grading.

Denote by $\CC[J_\infty\SSS](n)\subset\CC[J_\infty\SSS]$  the component corresponding to the character $t\mapsto t^n$. Then one has, by the definition of the
asymptotic vertex algebra,
\[
\cW^k(\fg,f)\stackrel{\sim}{\longrightarrow}\prod_{n\ll+\infty}\CC[J_\infty\SSS](n)\hbar^{-n/2},
\]
where the symbol $\prod_{n\gg-\infty}$ means only those infinite sums $\sum a_n\hbar^{-n/2}$, $a_n\in\CC[J_\infty\SSS](n)$ are allowed where $a_n=0$ for all $n$ sufficiently large. Now notice that since $\SSS=f+\fg^e$, sect.~\ref{inside-g},  elementary representation theory of $sl_2$ implies all weights of $\rho$ are nonnegative integers; hence
$\CC[J_\infty\SSS](n)=0$ if $n<0$. Therefore the above isomorphism becomes
\[
\cW^k(\fg,f)\stackrel{\sim}{\longrightarrow}\bigoplus_{n=0}^{+\infty}\CC[J_\infty\SSS](n)\hbar^{-n/2}.
\]
In other words, extracting $\CC^*$-invariants precludes the appearance of infinite sum, but this is what we had to prove. The case of $\cW^{-h^\vee}_0(\fg,f)=\cW^{-h^\vee}_0(\fg,f)_{\sqrt{\hbar}}^{\CC^*}$ is done in exactly the same way. \hfill $\qed$

\subsubsection{ }
\label{kac-roan-wak-version}
Kac, Roan, and Wakimoto define affine $W$-algebras attached to a not
necessarily even nilpotent element $f$ differently. 
Their way of dealing with the fact that $\chi=(f,.)$ is not a character of $\fg_{\geq 1}$ is to adjoin  an appropriate $\beta\gamma$-system rather
than restrict to a smaller subalgebra $\fm\subset \fg_{\geq 1}$. Let us recall their definition  \cite{KRW,KW}:

\begin{align*}
 \cW^k(\fg,f)_{KRW}=
H^{\infty/2+0}_{Q}(V^k(\fg){\otimes}
\cD_{\CC^r}^{ch}(\CC^r) {\otimes}
 Cl^{vert}(\fg_{\geq 1}\oplus 
\fg_{\geq 1}^*)),
\end{align*}
where
$r=\dim \fg_1/2$,
and 
\begin{align*}
 Q
= \sum_i ( u_i-\Phi_i){\otimes} \phi_i^*-
\frac{1}{2}\sum_{i,j,k}1{\otimes} c_{ij}^k(\phi_k)_{(-1)}((\phi^*_i)_{(-1)}(\phi^*_j)).
\end{align*}
Here, $\fg_{\geq m}=\fg_m\oplus\fg_{m+1}\oplus\cdots$, cf. sect.\ref{inside-g} ,
 $\{u_i\}$ is a graded basis of $\fg_{\geq 1}$,
$\{c_{ij}^k\}$ are the structure constants relative to the chosen basis,
$ \Phi_i=
	 (f|u_i)$ for $u_i\in \fg_{2}$,
$ \Phi_i=
	 0$ for $u_i\in \fg_{\geq 3}$,
and $\{\Phi_i; u_i\in \fg_1\}$ is a set of strong generators of 
$\cD_{\CC^r}(\CC^r)$ of the form $x_\bullet$ or $\partial_{\bullet}$ s.t. $\Phi_{i(0)}\Phi_j=(f,[u_i,u_j])$. Note that the symplectic form
$u_i,u_j\mapsto(f,[u_i,u_j])$ has appeared in sect.\ref{inside-g} .

Likewise, cf.\cite{A3}, Proposition~2.2,
\begin{align*}
 \cW^{-h^\vee}_0(\fg,f)_{KRW}=
H^{\infty/2+0}_{Q}(V^{-h^\vee}_0(\fg){\otimes}
\cD_{\CC^r}^{ch}(\CC^r) {\otimes}
 Cl^{vert}(\fg_{\geq 1}\oplus 
\fg_{\geq 1}^*)),
\end{align*}

\subsubsection{ }
We wish to show the equivalence of the 2 definitions.  Since the definition of sect.~\ref{c-star-action-on-asympt-w and-affine-w} is what we actually need and what
arises naturally in the context of the present work, we will be somewhat
brief. Note that
it is not obvious whether  
the definition of 
of sect.~\ref{c-star-action-on-asympt-w and-affine-w} 
gives a quantization of
$J_{\infty}\mathbb{S}$ 
in a non-$\hbar$-adic setting,
because argument of sec.\ \ref{proofofthmchquantham-i-ii-iii}
does not work without the $\hbar$-adic completeness.
Note also that the equivalence of the $2$ definitions proves
that
the definition of $W$-algebras given in \cite{KRW}
is independent of the choice of the good grading,
which is a fact remarked  in \cite[Remark 5.4]{DeSole-Kac}
without proof.

 To begin with, we will introduce a new asymptotic vertex algebra that will serve as an intermediary. Extensively used above is the Clifford vertex algebra $Cl^{vert}(V\oplus V^*)$.
Likewise, given 2 vector spaces $V$ and $W$ with pairing $\langle.,.\rangle: V\times W\rightarrow \CC$ one defines in an obvious manner a vertex algebra $Cl^{vert}(V,W,\langle.,.\rangle)$, a ``skewed'' Clifford algebra.

We will be concerned with this construction when $V=\fg_{\geq 1}^*$, $W=\fg_{\geq 2}$, and $\langle.,.\rangle$ is the natural pairing of $\fg_{\geq 1}^*$ and $\fg_{\geq 1}$
restricted to $\fg_{\geq 2}\subset \fg_{\geq 1}$. Denote the corresponding skewed Clifford algebra simply by $Cl^{vert}(\fg_{\geq 1}^*,\fg_{\geq 2})$. 

\begin{sloppypar}
Now consider  $V^k(\fg)_{\sqrt{\hbar}}\hat{\otimes} Cl^{vert}(\fg_{\geq 1}^*\otimes\fg_{\geq 1})_\hbar$. It contains a subalgebra,
$V^k(\fg)_{\sqrt{\hbar}}\hat{\otimes} Cl^{vert}(\fg_{\geq 1}^*,\fg_{\geq 2})_\hbar$, and a familiar element
\end{sloppypar}
\[
\tilde{d}^{ch}_{\chi}= \sum_i\hbar^{-1}( g_i-(f,g_i))\otimes \phi_i^*-
\frac{1}{2}\sum_{i,j,k}1\otimes c_{ij}^k(\phi_k)_{(-1)}((\phi^*_i)_{(-1)}(\phi^*_j)).
\]
It differs from its twin $d^{ch}_{\chi}$ in (\ref{chir-ham-red-compl-usu-defn}) only in that it is based on a bigger Lie algebra, $\fg_{\geq 1}$; in particular, $\{g_i\}$ is a basis
of the latter. Of course $\chi=(f,.)\in\fg^*_{\geq 1}$, but it is not a character and therefore $(\tilde{d}^{ch}_{\chi})_{(0)}^2\neq 0$. Nevertheless, it is a derivation (any $v_{(0)}$ is),
it preserves the subalgebra $V^k(\fg)_{\sqrt{\hbar}}\hat{\otimes}
Cl^{vert}(\fg_{\geq 1}^*,\fg_{\geq 2})_\hbar$, and its restriction to
the latter does square to 0
 because of the fact that $\chi$ is a character of $\fg_{\geq 2}$.
 Thus we obtain
a differential graded asymptotic vertex algebra $(V^k(\fg)_{\sqrt{\hbar}}\hat{\otimes} Cl^{vert}(\fg_{\geq 1}^*,\fg_{\geq 2})_\hbar,(\tilde{d}^{ch}_{\chi})_{(0)})$ along with its
``reduced'' version $(V^{-h^\vee}_0(\fg)_{\sqrt{\hbar}}\hat{\otimes} Cl^{vert}(\fg_{\geq 1}^*,\fg_{\geq 2})_\hbar,(\tilde{d}^{ch}_{\chi})_{(0)})$.

\subsubsection{ }
\label{hadic-krw-and-thm-on-equi}
The Kac-Roan-Wakimoto construction has an obvious asymptotic version: 
\[(V^k(\fg)_{\sqrt{\hbar}}\hat{\otimes} \cD_{\CC^r}^{ch}(\CC^r)_\hbar {\otimes} Cl^{vert}(\fg_{\geq 1}\oplus \fg_{\geq 1}^*)_\hbar, Q_\hbar)\text{ or } 
(V^{-h^\vee}_0(\fg)_{\sqrt{\hbar}}\hat{\otimes} \cD_{\CC^r}^{ch}(\CC^r)_\hbar {\otimes} Cl^{vert}(\fg_{\geq 1}\oplus \fg_{\geq 1}^*)_\hbar, Q_\hbar)
\]
 that lead to
 \[
 \cW^k(\fg,f)_{KRW,\sqrt{\hbar}}=H^{\infty/2+0}_{Q_\hbar}(V^k(\fg)_{\sqrt{\hbar}}\hat{\otimes} \cD_{\CC^r}^{ch}(\CC^r)_\hbar {\otimes} Cl^{vert}(\fg_{\geq 1}\oplus \fg_{\geq 1}^*)_\hbar),
 \]
 \[
 \cW^{-h^\vee}_0(\fg,f)_{KRW,\sqrt{\hbar}}=H^{\infty/2+0}_{Q_\hbar}(V^{-h^\vee}_0(\fg)_{\sqrt{\hbar}}\hat{\otimes} \cD_{\CC^r}^{ch}(\CC^r)_\hbar {\otimes} Cl^{vert}(\fg_{\geq 1}\oplus \fg_{\geq 1}^*)_\hbar).
 \]
The indicated differential graded  asymptotic vertex algebras carry a $\CC^*$ action (defined as in sect.~\ref{c-star-action-on-asympt-w and-affine-w}) and the canonical filtration
$F$ intoduced in (\ref{filtr-on-acdo}). One has isomorphisms
\begin{equation}
\label{prop-krw-hadic-1}
 \cW^k(\fg,f)_{KRW,\sqrt{\hbar}}^{\CC^*}\stackrel{\sim}{\longrightarrow} \cW^k(\fg,f)_{KRW},\; 
  \cW^{-h^\vee}_0(\fg,f)_{KRW,\sqrt{\hbar}}^{\CC^*}\stackrel{\sim}{\longrightarrow}  \cW^{-h^\vee}_0(\fg,f)_{KRW} ,
  \end{equation}
 \begin{equation}
\label{prop-krw-hadic-2}
Gr_{F}H^{\infty/2+\bullet}_{Q_\hbar}(V^k(\fg)_{\sqrt{\hbar}}\hat{\otimes} \cD_{\CC^r}^{ch}(\CC^r)_\hbar {\otimes} Cl^{vert}(\fg_{\geq 1}\oplus \fg_{\geq 1}^*)_\hbar)
\stackrel{\sim}{\longrightarrow}\CC[J_\infty\SSS][\hbar,\hbar^{-1}],
\end{equation}
\begin{equation}
\label{prop-krw-hadic-3}
Gr_{F}H^{\infty/2+\bullet}_{Q_\hbar}(V^{-h^\vee}_0(\fg)_{\sqrt{\hbar}}\hat{\otimes} \cD_{\CC^r}^{ch}(\CC^r)_\hbar {\otimes} Cl^{vert}(\fg_{\geq 1}\oplus \fg_{\geq 1}^*)_\hbar)
\stackrel{\sim}{\longrightarrow} \CC[J_\infty S][\hbar,\hbar^{-1}].
  \end{equation}
 The two latter isomorphisms follow easily from \cite{A3}, e.g., Proposition~2.3 in {loc. cit.}

We have a triple of subalgebras of $\fg$: $\fg_{\geq 2}\subset\fm\subset\fg_{\geq 1} $.  The natural maps $\fg_{\geq 1}^*\rightarrow \fm^*$ and $\fg_{\geq 2}\hookrightarrow\fm$
induce asymptotic vertex algebra morphisms
\begin{equation}
\label{interm-we-1}
V^k(\fg)_{\sqrt{\hbar}}\hat{\otimes} Cl^{vert}(\fg_{\geq 1}^*,\fg_{\geq 2})_\hbar\rightarrow V^k(\fg)_{\sqrt{\hbar}}\hat{\otimes} Cl^{vert}(\fm^*\oplus\fm)_\hbar,
\end{equation}
\begin{equation}
\label{interm-we-2}
V^{-h^\vee}_0(\fg)_{\sqrt{\hbar}}\hat{\otimes} Cl^{vert}(\fg_{\geq 1}^*,\fg_{\geq 2})_\hbar\rightarrow V^{-h^\vee}_0(\fg)_{\sqrt{\hbar}}\hat{\otimes} Cl^{vert}(\fm^*\oplus\fm)_\hbar ,
\end{equation}
both of which are readily seen to respect the differentials.

Similarly, the embedding $\fm\hookrightarrow \fg_{\geq 1}$ induces  differential graded asymptotic vertex algebra morphisms
\begin{equation}
\label{interm-krw-1}
V^k(\fg)_{\sqrt{\hbar}}\hat{\otimes} Cl^{vert}(\fg_{\geq 1}^*,\fg_{\geq 2})_\hbar\rightarrow V^k(\fg)_{\sqrt{\hbar}}\hat{\otimes} \cD_{\CC^r}^{ch}(\CC^r)_\hbar \hat{\otimes} Cl^{vert}(\fg_{\geq 1}\oplus \fg_{\geq 1}^*)_\hbar,
\end{equation}
\begin{equation}
\label{interm-krw-2}
V^{-h^\vee}_0(\fg)_{\sqrt{\hbar}}\hat{\otimes} Cl^{vert}(\fg_{\geq 1}^*,\fg_{\geq 2})_\hbar\rightarrow V^{-h^\vee}_0(\fg)_{\sqrt{\hbar}}\hat{\otimes} \cD_{\CC^r}^{ch}(\CC^r)_\hbar \hat{\otimes} Cl^{vert}(\fg_{\geq 1}\oplus \fg_{\geq 1}^*)_\hbar.
\end{equation}
\begin{thm}
\label{thm-on-equi-of-2-defs}
The 4 morphisms (\ref{interm-we-1}--\ref{interm-krw-2}) are quasiisomorphisms.
\end{thm}
Thanks to (\ref{prop-krw-hadic-1}) the following is immediate:
\begin{cor}
\label{cor-on-equi-of-2-defs}
\[
\cW^k(\fg,f)_{\sqrt{\hbar}}\stackrel{\sim}\longrightarrow \cW^k(\fg,f)_{KRW,\sqrt{\hbar}},\;
\cW^{-h^\vee}_0(\fg,f)_{\sqrt{\hbar}}\stackrel{\sim}\longrightarrow \cW^{-h^\vee}_0(\fg,f)_{KRW,\sqrt{\hbar}},
\]
\[
\cW^k(\fg,f)\stackrel{\sim}\longrightarrow \cW^k(\fg,f)_{KRW},\;
\cW^{-h^\vee}_0(\fg,f)\stackrel{\sim}\longrightarrow \cW^{-h^\vee}_0(\fg,f)_{KRW}.
\]
\end{cor}  

{\em Proof of Theorem~\ref{thm-on-equi-of-2-defs}.} Begin with morphisms (\ref{interm-we-1}) and (\ref{interm-we-2}) .  Consider the induced morphism of
graded objects
\begin{equation}
\label{what-to-prove-equi-def}
Gr_F H^{\infty/2+\bullet}_{\tilde{d}^{ch}_\chi}(V^k(\fg)_{\sqrt{\hbar}}\hat{\otimes} Cl^{vert}(\fg_{\geq 1}^*,\fg_{\geq 2})_\hbar)
\rightarrow 
Gr_F H^{\infty/2+\bullet}_{d^{ch}_\chi}(V^k(\fg)_{\sqrt{\hbar}}\hat{\otimes} Cl^{vert}(\fm^*\oplus\fm)_\hbar).
\end{equation}
Converging to each of these is the spectral sequence determined by the $\hbar$-adic filtration of the original complexes, 
the one utilized in the proof of Theorem~\ref{chiral-quant-hamilt-red}.  Furthermore, morphism (\ref{interm-we-1})  determines a morphism of spectral sequences.
The spectral sequence for the r.h.s of (\ref{what-to-prove-equi-def}) has been dealt with in Lemma~\ref{why-quantizs-of-jet-slice}: it collapses and  gives
$\CC[J_\infty\SSS][\hbar,\hbar^{-1}]$ sitting in the cohomological degree 0. The computation of the spectral sequence for the l.h.s of (\ref{what-to-prove-equi-def})  can be derived in the same way
as that for the r.h.s. from Theorem~\ref{chiral-quant-hamilt-red} except that while in the latter case one uses the Gan-Ginzburg isomorphism (\ref{GG-Lemma}), the former is based
on a similar and easy to verify isomorphism
\[
M_{\geq 1}\times\SSS\stackrel{\sim}{\longrightarrow} f+\fg_{\geq 2}^\perp.
\]
The answer is again $\CC[J_\infty\SSS][\hbar,\hbar^{-1}]$ sitting in the cohomological degree 0. Therefore (\ref{what-to-prove-equi-def}) is an isomorphism and so is 
(\ref{interm-we-1}).
 Morphism  (\ref{interm-we-2}) is handled in exactly the same way.
 
 The argument on the Kac-Roan-Wakimoto side,(\ref{interm-krw-1}) and (\ref{interm-krw-2}), is similar. One passes to the induced map of graded objects, and when
 dealing with the r.h.s. instead of Theorem~\ref{chiral-quant-hamilt-red}  one relies on  (\ref{prop-krw-hadic-2}) and (\ref{prop-krw-hadic-3}).\hfill $\qed$

By \cite{A3},
$\cW_0^{-h^{\vee}}(\fg,f)_{KRW}$ is simple. 
Therefore, we have the following
\begin{cor}
The vertex algebra
$\cW_0^{-h^{\vee}}(\fg,f)$ is simple.
\end{cor}

\subsubsection{ }
\label{acdo-on-S-intro}
To conclude our discussion of various (asymptotic) vertex algebras involved,
let us quantize $\Gamma(J_\infty\tilde{S}, \cO_{J_\infty\tilde{S}})$. It is known that $ch_2(X)=0$, therefore, Theorem~\ref{thn-on-class-cdo},
$X$ carries a CDO, which is unique up to isomorphism \cite{GMSII};
denote it by $\cD^{ch}_X$.  We obtain then, sect.~\ref{main-def}, an ACDO $\cD^{ch}_{X,\hbar}$ with lattice $\cD^{ch}_{X, \hbar,\geq}$, cf.
Definition~\ref{defn-of-hafic-asympt-cdo}, and the chiral Hamiltonian reduction complex
\begin{eqnarray}
\label{chir-ham-red-compl-acdo-defn}
&(\cD^{ch}_{X,\hbar}\hat{\otimes} Cl^{vert}(\fm\oplus\fm^*)_\hbar, (d^{ch}_{\chi})_{(0)}),\\
& d^{ch}_{\chi}= \sum_i( \mu_{ch}(m_i)-\frac{(f,m_i)}{\hbar})\otimes \phi_i^*-
\frac{1}{2}\sum_{i,j,k}1\otimes c_{ij}^k(\phi_k)_{(-1)}((\phi^*_i)_{(-1)}(\phi^*_j)).\nonumber
\end{eqnarray}
By Theorem~\ref{chiral-quant-hamilt-red}, the sheaf
\begin{equation}
\label{chir-ham-red-for-preim-slow-slice}
\cD^{ch}_{\tilde{S},\hbar}\stackrel{\text{def}}{=}p_{\infty*}\cH^{\infty/2+0}_{d^{ch}_\chi}(\cD^{ch}_{X,\hbar}\hat{\otimes} Cl^{vert}(\fm\oplus\fm^*)_\hbar),\; \chi=(f,.)
\end{equation}
is an ACDO and a quantization of $\cO_{J_\infty\tilde{S}}$. We wish to  relate $\Gamma(J_\infty\tilde{S},\cD^{ch}_{\tilde{S},\hbar})$ and
$\cW^{-h^\vee}_0(\fg,f)$.

\subsection{Main result}
\label{main-result}
\subsubsection{ }
\label{c-star-equi}
The CDO $\cD^{ch}_X$ is known to carry an affine Lie algebra action at the critical level, i.e., there is a morphism
$V_0^{-h^\vee}(\fg)\rightarrow \Gamma(X,\cD^{ch}_{X})$, \cite{MSV,GMSII}. Hence, Corollary~\ref{from-emd-in-cdo-to-emb-in-acdo}, a morphism
\begin{equation}
\label{aff-alg-str-on-acdo}
\phi: V_0^{-h^\vee}(\fg)\rightarrow \Gamma(J_\infty T^*X,\cD^{ch}_{X,\hbar}).
\end{equation}
Note that its restriction to $V(\fm)\subset V^{-h^{\vee}}(\fg)$ coincides with $\mu^{ch}$, by definition.
 In particular,(\ref {aff-alg-str-on-acdo}) gives us a Lie algebra morphism
\[
sl_2\rightarrow\text{End}(\cD^{ch}_{X,\hbar}),\;
\left(\begin{array}{rr}a&b\\c&-a\end{array}\right)\mapsto a\phi(h)_{(0)}+b\phi(e)_{(0)}+c\phi(f)_{(0)},
\]
where $e$, $h$, $f$ is the $sl_2$-triple introduced in sect.~\ref{inside-g}.

By construction, this action is a lift of the canonical action on $X$. Therefore, it can be integrated to give $\cD^{ch}_{X,\hbar}$ an $SL_2$-equivariant
sheaf structure. Let
\[
SL_2\ni g\mapsto \gamma(g): \cD^{ch}_{X,\hbar}(g^{-1}U)\rightarrow\cD^{ch}_{X,\hbar}(U),\; \forall U\subset J_\infty T^*X
\]
be the map defining this $SL_2$-equivariant structure. We obtain then  a $\CC^*$-equivariant structure by restricting  as follows:
\[
\CC^*\ni t\mapsto \rho_1(t)=\gamma(\left(\begin{array}{ll}t^{-1}&0\\0&t\end{array}\right)).
\]
The ACDO $\cD^{ch}_{X,\hbar}$ carries yet another $\CC^*$-equivariant structure, one induced by the square of the canonical $\CC^*$-action on
the vector bundle $T^*X\rightarrow X$, meaning $t\cdot\xi=t^2\xi$ for $t\in\CC^*$, $\xi\in\cT_X\subset\cO_{J_\infty T^*X}$. We assume given
an atlas of $X$, where locally  $\cD^{ch}_{X,\hbar}$ is generated by functions and `coordinate vector fields' $\{\widehat{\partial}_i\}$ subject to
standard relations, see the proof of Lemma~\ref{quasiclass-lim-hadic-cdo} for explicit formulas. Now define 
\[
t\cdot\widehat{\partial}_i=t^2\widehat{\partial}_i, t\cdot\hbar=t^2\hbar.
\]
The formulas recorded in the proof of Lemma~\ref{quasiclass-lim-hadic-cdo} show at once that this action respects the relations and the
transformation functions and, therefore, defines 
\[
\CC^*\ni t\mapsto \rho_2(t): \cD^{ch}_{X,\hbar}(t^{-1}\cdot U)\rightarrow\cD^{ch}_{X,\hbar}(U),\; \forall U\subset J_\infty T^*X.
\]

It is clear that $\rho_1(t_1)\circ\rho_2(t_2)=\rho_2(t_2)\circ\rho_1(t_1)$, and we define the desired $\CC^*$-action
 \[
 \CC^*\ni t\mapsto \rho(t)\stackrel{\text{def}}{=}\rho_1(t)\circ\rho_2(t).
 \]
 All of this is comletely analogous to sect.~\ref{c-star-action-on-asympt-w and-affine-w}.

Note again that $\rho(t)(\hbar)=t^2\hbar$.
It is convenient to have $\hbar$ generate the identity character $\{t\mapsto t\}$ of $\CC^*$.  To achieve this, extend the scalars and define
\begin{equation}
\label{defn-acdo-after-ext-scal}
\cD^{ch}_{X,\sqrt{\hbar}}=\cD^{ch}_{X,\hbar}\otimes_{\CC((\hbar))}\CC((\sqrt{\hbar}))
\end{equation}
 stipulating that $\rho(t)(\sqrt{\hbar})=t\sqrt{\hbar}$.
This notation is clearly ambiguous, but we will use it for the lack of a better one.

\subsubsection{ }
\label{equi-strre-on-the-complex}
Consider the chiral Hamiltonian reduction complex $(\cD^{ch}_{X,\sqrt{\hbar}}\hat{\otimes} Cl^{vert}(\fm\oplus\fm^*)_\hbar, (d^{ch}_{\chi})_{(0)})$,  an obvious extension of
(\ref{chir-ham-red-compl-acdo-defn}) with the same differential. Define the $\CC^*$-equivariant structure
\[
\CC^*\ni t\mapsto \rho(t)\otimes t^{-h}: \cD^{ch}_{X,\sqrt{\hbar}}(\rho(t^{-1})U)\hat{\otimes} Cl^{vert}(\fm\oplus\fm^*)_\hbar\rightarrow\cD^{ch}_{X,\sqrt{\hbar}}(U)\hat{\otimes}
 Cl^{vert}(\fm\oplus\fm^*)_\hbar,
\]
where $h$ acts naturally on $Cl^{vert}(\fm\oplus\fm^*)$.
This makes $\Gamma(J_\infty T^*X,\cD^{ch}_{X,\sqrt{\hbar}}\hat{\otimes} Cl^{vert}(\fm\oplus\fm^*)_\hbar)$ a $\CC^*$-module, and it is easy to see that the differential
$d^{ch}_{\chi}$ is $\CC^*$-invariant. (Note that a more obvious action $\rho\otimes 1$ does not do the job.) Hence the complex
$(\cD^{ch}_{X,\sqrt{\hbar}}\hat{\otimes} Cl^{vert}(\fm\oplus\fm^*)_\hbar, (d^{ch}_{\chi})_{(0)})$ is $\CC^*$-equivariant.  Thus we obtain a $\CC^*$-equivariant structure
on the cohomology $\cH^{\infty/2+\bullet}_{d^{ch}_{\chi}}(\cD^{ch}_{X,\sqrt{\hbar}}\hat{\otimes} Cl^{vert}(\fm\oplus\fm^*)_\hbar)$ and on 
\begin{equation}
\label{vers-acdo-tilde-hsqrt}
\cD^{ch}_{\tilde{S},\sqrt{\hbar}}\stackrel{\text{def}}{=}p_{\infty *}\cH^{\infty/2+0}_{d^{ch}_{\chi}}(\cD^{ch}_{X,\sqrt{\hbar}}\hat{\otimes} Cl^{vert}(\fm\oplus\fm^*)_\hbar),
\end{equation}
cf. (\ref{chir-ham-red-for-preim-slow-slice}).

\subsubsection{ }
\label{we-form-main-res}
Morphism (\ref{aff-alg-str-on-acdo})  extends to a morphism of the corresponding asymptotic vertex algebras
\[
\phi_\hbar: V_0^{-h^\vee}(\fg)_{\sqrt{\hbar}} \rightarrow \Gamma(J_\infty T^*X,\cD^{ch}_{X,\sqrt{\hbar}}).
\]
and defines an obvious morphism of chiral Hamiltonian reduction complexes  (\ref{chir-ham-red-compl-usu-defn}) and (\ref{chir-ham-red-compl-acdo-defn})
\[
 V_0^{-h^{\vee}}(\fg)_{\sqrt{\hbar}}\hat{\otimes} Cl^{vert}(\fm\oplus\fm^*)_\hbar\rightarrow \Gamma(J_\infty T^*X,\cD^{ch}_{X,\sqrt{\hbar}}\hat{\otimes} Cl^{vert}(\fm\oplus\fm^*)_\hbar).
\]
Furthermore, an obvious comparison of sect.~\ref{c-star-action-on-asympt-w and-affine-w} and \ref{c-star-equi} shows that it respects the $\CC^*$-module structures.

The passage to the cohomology defines
\[
\cW_0^{-h^{\vee}}(\fg,f)_{\sqrt{\hbar}}\rightarrow H^{\infty/2+0}_{d^{ch}_\chi}(\Gamma(J_\infty T^*X,\cD^{ch}_{X,\sqrt{\hbar}}\hat{\otimes} Cl^{vert}(\fm\oplus\fm^*)_\hbar)),
\]
\[
\cW_0^{-h^{\vee}}(\fg,f)\rightarrow H^{\infty/2+0}_{d^{ch}_\chi}(\Gamma(J_\infty T^*X,\cD^{ch}_{X,\sqrt{\hbar}}\hat{\otimes} Cl^{vert}(\fm\oplus\fm^*)_\hbar))^{\CC^*}.
\]

The definition of $\cD^{ch}_{\tilde{S},\hbar}$ as the sheaf associated to the presheaf of chiral Hamiltonian reductions, sect.~\ref{quant-chital-hamilt-red},
implies a morphism
\[
H^{\infty/2+0}_{d^{ch}_\chi}(\Gamma(J_\infty T^*X,\cD^{ch}_{X,\sqrt{\hbar}}\hat{\otimes} Cl^{vert}(\fm\oplus\fm^*)_\hbar))\rightarrow
\Gamma(J_\infty \tilde{S},\cD^{ch}_{\tilde{S},\sqrt{\hbar}}).
\]
Composing with the previous ones we obtain
\begin{equation}
\label{morph-source-main-thm}
\cW^{-h^\vee}_0(\fg,f)_{\sqrt{\hbar}}\rightarrow\Gamma(J_{\infty}\tilde{S},\cD^{ch}_{\tilde{S},\sqrt{\hbar}}),\;
\cW^{-h^\vee}_0(\fg,f)\rightarrow\Gamma(J_{\infty} \tilde{S},\cD^{ch}_{\tilde{S},\sqrt{\hbar}})^{\CC^*}.
\end{equation}
Here is the main result of this note:
\begin{thm}
\label{main-res}
Morphisms (\ref{morph-source-main-thm}) is an isomorphism.
\end{thm}

\subsubsection{Proof.}
\label{proof-of-main-res}

Consider the natural composite projection $\pi: J_\infty T^*X\rightarrow T^*X\rightarrow X$ and $\{U_i\}$, an open affine cover of $X$. Denote
by $\fU$ the affine cover $\{\pi^{-1} U_i\}$. For any sheaf  $\cA$ on $J_\infty X$, denote by $C^\bullet(\fU,\cA)$ the corresponding \v Cech complex. 

For any $V(\fm)_\hbar$-module $M$, let $C^{\infty/2+\bullet}(\fm, \chi;M)$ stand for the chiral Hamiltonian reduction complex
$(M\hat{\otimes} Cl^{vert}(\fm\oplus\fm^*)_\hbar, (d^{ch}_\chi)_{(0)})$; the familiar complex $C^{\infty/2+\bullet}(\fm, \chi;\cD^{ch}_{X,\sqrt{\hbar}})$ will be of special
interest for us. 

\begin{sloppypar}
Thus we obtain a bicomplex  
$C^\bullet(\fU,C^{\infty/2+\bullet}(\fm, \chi;\cD^{ch}_{X,\sqrt{\hbar}}))$.
Denote by  $H^\bullet_{tot}(\fU,C^{\infty/2+\bullet}(\fm,\chi; \cD^{ch}_{X,\sqrt{\hbar}}))$ its cohomology w.r.t. the total differential. Define the hypercohomology as usual
\end{sloppypar}
\[
\HH^{\infty/2+\bullet}(J_\infty T^*X, C^{\infty/2+\bullet}(\fm, \chi;\cD^{ch}_{X,\sqrt{\hbar}}))=\lim_{\fU}
H_{tot}^\bullet(\fU,C^{\infty/2+\bullet}(\fm, \chi;\cD^{ch}_{X,\sqrt{\hbar}})).
\]
\begin{lem}
\label{from-hypercoho-to-main-res}
There exist isomorphisms
\[
\Gamma(J_\infty \tilde{S},\cD^{ch}_{\tilde{S},\sqrt{\hbar}})\stackrel{\sim}{\longleftarrow}\HH^{\infty/2+0}(J_\infty T^*X, C^{\infty/2+\bullet}(\fm,\chi;
 \cD^{ch}_{X,\sqrt{\hbar}}))
\stackrel{\sim}{\longrightarrow}\cW^{-h^\vee}_0(\fg,f).
\]
\end{lem}
\begin{sloppypar}
{\em Proof of Lemma~\ref{from-hypercoho-to-main-res}.}  
For each $\fU$, there are 2 spectral sequences, 
$\{`E^{pq}_r\}$ and $\{``E^{pq}_r\}$, the former
starting with $(d^{ch}_\chi)_{(0)}$, the latter with $d_{\check{C}}$. Both  converge to $H_{tot}^\bullet(\fU,C^{\infty/2+\bullet}(\fm, \chi;\cD^{ch}_{X,\sqrt{\hbar}}))$ because
the complex is of finite length in the \v Cech direction. Theorem~\ref{chiral-quant-hamilt-red}  implies that  (i) $`E^{pq}_1=0$ if $q<0$ (we are assuming that $q$ is the
coordinate in the Hamiltonian reduction complex direction), and (ii)
 $`E^{00}_2=\Gamma(J_\infty
 T^*X,\cH^{\infty/2+0}_{d^{ch}_\chi}(\cD^{ch}_{X,\hbar}\hat{\otimes}
 Cl^{vert}(\fm\oplus\fm^*)_\hbar))$.  
Furthermore,  by the definition of the direct image, $\Gamma(J_\infty T^*X,\cH^{\infty/2+0}_{d^{ch}_\chi}(\cD^{ch}_{X,\hbar}\hat{\otimes} Cl^{vert}(\fm\oplus\fm^*)_\hbar))=\Gamma(J_{\infty}\tilde{S},(p_\infty)_*\cH^{\infty/2+0}_{d^{ch}_\chi}(\cD^{ch}_{X,\hbar}\hat{\otimes} Cl^{vert}(\fm\oplus\fm^*)_\hbar))$, and so, by virtue of
(\ref{chir-ham-red-for-preim-slow-slice}), (ii) becomes 
$`E^{00}_2=\Gamma(J_{\infty}\tilde{S},\cD^{ch}_{\tilde{S},\hbar})$

Assertion (i) implies that
no higher differentials touch $`E^{00}_2$ and no terms different from $`E^{00}_2$ contribute to $\HH^{\infty/2+0}(J_\infty T^*X, C^{\infty/2+\bullet}(\fm,\chi;
 \cD^{ch}_{X,\sqrt{\hbar}}))$, hence the leftmost isomorphism of the lemma.
\end{sloppypar}

Now focus on $\{``E^{pq}_r\}$. We have $``E_1^{pq}= C^{\infty/2+q}(\fm,\chi;H^p(J_\infty T^*X, \cD^{ch}_{X,\sqrt{\hbar}}))$.  It was proved in \cite{AChM}
that $H^0(X, \cD^{ch}_X)=V^{-h^\vee}_0(\fg)$. In the present situation literally the same proof gives
$
H^0(J_{\infty}T^*X, \cD^{ch}_{X,\sqrt{\hbar}})=V^{-h^\vee}_0(\fg)_{\sqrt{\hbar}}.
$
Therefore, by Lemma~\ref{why-quantizs-of-jet-slice}, $``E_2^{0\bullet}= \cW^{-h^\vee}_0(\fg,f)$  sitting in cohomological degree 0.  We would like to show
that the natural map $``E_2^{0\bullet}\rightarrow \HH^{\infty/2+0}(J_\infty T^*X, C^{\infty/2+\bullet}(\fm,\chi;
 \cD^{ch}_{X,\sqrt{\hbar}}))$ is an isomorphism.  In order to do so, we need to get hold of $H^p(J_\infty T^*X, \cD^{ch}_{X,\sqrt{\hbar}})$. A similar cohomology,
$H^p(X, \cD^{ch}_{X})$, was computed in  \cite{AM}, and we have to analyze the effect of introducing  $\hbar$.

In order to unburden the notation write $C^\bullet$ for $C^\bullet(\fU,C^{\infty/2+\bullet}(\fm, \chi;\cD^{ch}_{X,\sqrt{\hbar}}))$. Filtration (\ref{filtr-on-acdo}) on
$C^{\infty/2+\bullet}(\fm, \chi;\cD^{ch}_{X,\sqrt{\hbar}})$ defines a filtration $\{F_\bullet C^\bullet\}$ and an exact sequence
\[
0\longrightarrow F_NC^\bullet\longrightarrow C^\bullet\longrightarrow C^\bullet/F_NC^\bullet\longrightarrow 0.
\]
Hence an exact sequence of hypercohomology
\begin{equation}
\label{why-h-divisible-1}
\HH^0(F_NC^\bullet)\longrightarrow \HH^0(C^\bullet)\longrightarrow \HH^0(C^\bullet/F_NC^\bullet).
\end{equation}
Each of these hypercohomology groups can be computed by a version of either of the 2 spectral sequences introduced above for $C^\bullet$. Focus
on  the 2nd one: $``E^{pq}_r(F_NC^\bullet)$, $``E^{pq}_r(C^\bullet)$, which appeared above simply as $``E^{pq}_r$, and $``E^{pq}_r(C^\bullet/F_NC^\bullet)$.
The projection $C^\bullet\rightarrow C^\bullet/F_NC^\bullet$ induces a collection  of  sequences (compatible with differentials) 
\begin{equation}
\label{why-h-divisible-2}
``E^{pq}_r(F_NC^\bullet)\longrightarrow ``E^{pq}_r(C^\bullet)\longrightarrow``E^{pq}_r(C^\bullet/F_NC^\bullet).
\end{equation}
 According to \cite{AM},
$H^\bullet(X,\cD^{ch}_{X})= H^\bullet(\fn,\CC)\otimes V^{-h^\vee}_0(\fg)$.  Now one may  want to obtain an analogue: 
$H^\bullet(J_\infty T^*X,\cD^{ch}_{X,\sqrt{\hbar}})= H^\bullet(\fn,\CC)\otimes V^{-h^\vee}_0(\fg)_{\sqrt{\hbar}}$, or equivalently,
$``E^{pq}_1(C^\bullet)=C^{\infty/2+q}(\fm,\chi;\oplus V^{-h^\vee}_0(\fg)_{\sqrt{\hbar}})$. It is not quite clear whether this is correct
because of the infinite power series in $\hbar$ involved, but what is easy to see is the validity for the truncated version $C^\bullet/F_NC^\bullet$:
\[
``E^{pq}_1(C^\bullet/F_NC^\bullet)=C^{\infty/2+q}(\fm,\chi;\oplus V^{-h^\vee}_0(\fg)_{\sqrt{\hbar}})/F_NC^{\infty/2+q}(\fm,\chi;\oplus V^{-h^\vee}_0(\fg)_{\sqrt{\hbar}}).
\]
Now the vanishing result of Lemma~\ref{why-quantizs-of-jet-slice}, or rather its proof,  implies that for each $N$
the spectral sequence $``E^{pq}_r(C^\bullet/F_NC^\bullet)$ collapses: $``E^{pq}_2(C^\bullet/F_NC^\bullet)=0$ if $q\neq 0$. Therefore  (\ref{why-h-divisible-1}) and
(\ref{why-h-divisible-2})
imply
\[
\oplus_{p\neq 0}``E^{p,-p}_\infty(C^\bullet)\subset \bigcap_{N}Gr F_N\HH^0(C^\bullet).
\]
But we already know what $\HH^0(C^\bullet)$ is: by the leftmost isomorphism of  Lemma~\ref{from-hypercoho-to-main-res} proved above,
$\HH^0(C^\bullet)=\Gamma(J_\infty \tilde{S},\cD^{ch}_{\tilde{S},\sqrt{\hbar}})$; it is, therefore, a subalgebra of an asymptotic vertex algebra, hence
$ \cap_N F_N\HH^0(C^\bullet)=0$. Hence $\oplus_{p\neq 0}``E^{p,-p}_\infty(C^\bullet)=0$  and
$\cW^{-h^\vee}_0(\fg,f)=``E_2^{0\bullet}\rightarrow \HH^{\infty/2+0}(J_\infty T^*X, C^{\infty/2+\bullet}(\fm,\chi;
 \cD^{ch}_{X,\sqrt{\hbar}}))$ is at least a surjection. But  $\cW^{-h^\vee}_0(\fg,f)$ is simple \cite{A3} hence the latter map is an isomorphism.
 \hfill $\qed$
 
 It is clear that the composition of the 2 isomorphisms of Lemma~\ref{from-hypercoho-to-main-res}  (one inverted) is precisely (\ref{morph-source-main-thm}).
 Theorem~\ref{main-res} has been proved. \hfill $\qed$
\begin{rem}
Notice that an analogue of (\ref{iso-of-funct-spaces})
\begin{equation}
\label{konst-thm-first-app-for-slodow-main-body}
\CC[J_\infty S]\longrightarrow\Gamma(J_\infty\tilde{S},\cO_{J_\infty\tilde{S}})
\end{equation}
is not an isomorphism, which makes the proof of Theorem~\ref{main-res} more involved but also vindicates the introduction of the quantum objects, such as
$\cD^{ch}_{\tilde{S},\sqrt{\hbar}}$.
\end{rem}

\subsection{Example: $\cW^{(2)}_3$}
\label{Example-of-sl3}

\newcommand{\NO}{{\genfrac{}{}{0pt}{1}{\circ}{\circ}}} 

As an illustration, we explicitly construct  in this section a localization of
the affine W-algebra $\cW^{(2)}_3$.
\subsubsection{ }
For $\fg=sl_3$ and $f = E_{21} \in sl_3$, 
the affine W-algebra $\cW^{k}(\fg, f)$ is known as
the Bershadsky-Polyakov algebra $\cW^{(2)}_3$,
see \cite{KRW}.
 The Bershadsky-Polyakov $\cW^{(2)}_3$ at level $k$ is a vertex algebra
 generated by the fields
$J(z)$, $G^{\pm}(z)$ and $S(z)$ with the following OPE's:
\begin{align}
\label{eq:ope-w23}
 &J(z)J(w)\sim \frac{2k +3}{3(z-w)^2},\quad G^{\pm}(z)G^{\pm}(w)\sim
 0,\\
&J(z)G^{\pm}(w)\sim \pm \frac{1}{z-w}G^{\pm}(w),\nonumber\\
&S(z)S(w)\sim
- \frac{(k+3)(2k+3)(3k+1)}{2(z-w)^4}+\frac{2(k+3)}{(z-w)^2}S(w)+\frac{k+3}{z-w}\partial
S(w),\nonumber\\
&S(z)G^{\pm}(w)\sim
 \frac{3(k+3)}{2(z-w)^2}G^{\pm}(w)+\frac{k+3}{z-w}\partial G^{\pm}(w),\nonumber\\
& S(z)J(w)\sim \frac{k+3}{(z-w)^2}J(w)+\frac{k+3}{z-w}\partial J(w),\nonumber\\
&G^{+}(z)G^-(w)\sim \frac{(k+1)(2k+3)}{(z-w)^3}
+\frac{3(k+1)}{(z-w)^2}J(w)\nonumber\\
&\qquad\qquad\qquad\qquad
+\frac{1}{z-w}\left(3 \NO J(w)^2 \NO + \frac{3(k+1)}{2}\partial J(w)
-S(w)\right). \nonumber
\end{align}
The vertex algebra $\cW^{(2)}_3$ was first introduced by
M.~Bershadsky \cite{Ber} and A.~Polyakov \cite{Pol}.

The direct calculation shows that the Feigin-Frenkel center of $\cW^{-h^{\vee}}(\fg, f)$
 is generated by $S(w)$ and the following field
\[
 S_3(z) = \NO G^-(z) G^+(z) \NO + \NO S(z) J(z) \NO 
 - \NO J(z)^3 \NO - 3 \NO J(z) \partial_z J(z) \NO - \partial_z^2 J(z).
\]

\subsubsection{} 
Let $\fg=sl_3$, $f=E_{21}\in sl_3$. We obtain an $sl_2$-triple $<f,h,e>$ with
$h=E_{11}-E_{22}$, $e=E_{12}$. The Slodowy slice is 
\[
\SSS\stackrel{\text{def}}{=}f+sl_3^e=\{\left(\begin{array}{ccc}
\delta&\alpha&\beta\\
1&\delta&0\\
0&\gamma&-2\delta
\end{array}
\right),\;\alpha,\beta,\gamma,\delta\in\CC\}.
\]
The Lie algebra $\fm$ can be chosen to be
\[
\fm=\{\left(\begin{array}{ccc}
0&u&v\\
0&0&0\\
0&0&0
\end{array}
\right),\;u,v\in\CC\}
\]
so that the corresponding group is 
\[
M=\{\left(\begin{array}{ccc}
1&u&v\\
0&1&0\\
0&0&1
\end{array}
\right),\;u,v\in\CC\}.
\]
The level set of the moment map, $(\mu^*)^{-1}(\chi)$, $\chi=(f,.)$, is
\[
(\mu^*)^{-1}(\chi)=
\{\left(\begin{array}{ccc}
h_1&x_{12}&x_{13}\\
1&h_2-h_1&x_{23}\\
0&x_{32}&-h_2
\end{array}
\right),\;h_\bullet, x_{\bullet\bullet}\in\CC\}.
\]
In this particular case a direct computation easily proves the Gan-Ginzburg result \cite{GG} that
for each element of $A\in (\mu^*)^{-1}(\chi)$ there is a unique element of $X \in M$ so that $XAX^{-1}\in\SSS$.
Therefore $\CC[\alpha,\beta,\gamma,\delta]$ is identified with $\CC[(\mu^*)^{-1}(\chi)/M]$, the Poisson reduced
variety. The same computation gives explicit formulas for $\alpha,\beta,\gamma,\delta$ as $M$-invariant functions on
$(\mu^*)^{-1}(\chi)$.  The answer is:
\begin{eqnarray}
\delta&=&\frac{1}{6}(E_{11}+E_{22}-2E_{33})\nonumber\\
\gamma&=&E_{32}\nonumber\\
\beta&=&E_{31}-(E_{11}-E_{33})E_{32}.\nonumber
\end{eqnarray}
(We let the curious reader  discover the formula for $\alpha$ on his own; we will not need this result.)

The intersection $S=\SSS\cap\cN$ is given by $P_2(X)=P_3(X)=0$, where $P_i(X)=\text{Tr}X^i$, $X\in\SSS$.
A quick computation shows that
\[
\CC[S]=\CC[\gamma,\beta,\delta]/\langle8\delta^3-\beta\gamma\rangle,
\]
which is a Kleinian singularity of type $A_2$, a particular case of a  famous theorem \cite{Bri,Slo}.

It is convenient to rescale by letting $h=2\delta$,  $a=\gamma$, $b=\beta$ so as to obtain
\[
\CC[S]=\CC[h,a,b]/\langle h^3-ab\rangle.
\]
The above formulas for $\gamma, \beta,\delta$   incidentally give (noncanonical) extensions of these functions to the entire $sl_3$, which
allows one to compute their Kirillov-Kostant Poisson brackets.  Upon subsequent restricting to $S$ one obtains the quadratic bracket
\[
\{h,a\}=a,\; \{h,b\}=-b,\; \{a,b\}=-3h^2;
\]
this is the Poisson reduced algebra.

Lemma~\ref{from-poiss-to-vert-poiss} defines the following vertex Poisson algebra structure on on the jet-algebra $\CC[J_\infty S]$:
\begin{equation}
\label{vert-poiss-on-jet-w23}
h_{(0)}a=a,\; h_{(0)}b=-b,\;a_{(0)}b=-3h^2,\; x_{(n)}y=0\text{ if }n>0,x,y=a,b\text{ or }h.
\end{equation}

\subsubsection{ }
\newcommand{\RR}{\mathbb R}

Consider the Springer resolution $\pi: \tilde{S} \longrightarrow S$.
It is well-known that the Slodowy variety $\tilde{S}$ coincides with
the minimal resolution of the Kleinian singularity $S$ (see \cite{Slo}).
Another realization of the minimal resolution is given as a quiver
variety which was studied in \cite{Kro}.
While a quantization of the Slodowy slice is given by the finite W-algebra,
a quantization of the quiver variety is called a deformed 
preprojective algebra or a rational Cherednik algebra.
Moreover, it is known that these two quantizations coincide as a result
of \cite{Pre} in the case of the Kleinian singularity of type $A$.
To construct an ACDO on $\tilde{S}$ explicitly,
we use a realization of the minimal resolution as a quiver variety.



For $i=1$, $2$, $3$, let $u_i$, $v_i$ be indeterminates. 
Set 
\[
 R = \CC[u_1, u_2, u_3, v_1, v_2, v_3] / 
\langle u_1 v_1 = u_2 v_2 = u_3 v_3 \rangle.
\]
The three-dimensional torus $T = \prod_{i \in \ZZ / 3 \ZZ} \CC^*$ acts on
$R$ by $g \cdot u_i = g_i^{-1} g_{i-1} u_i$, 
$g \cdot v_i = g_i g_{i-1}^{-1} v_i$ for $g = (g_j)_{j \in \ZZ / 3\ZZ}$
and $i=1$, $2$, $3$.
Namely $u_i$ is of $T$-weight $\varepsilon_{i-1} - \varepsilon_i$ and
$v_i$ is of $T$-weight $\varepsilon_i - \varepsilon_{i-1}$ where
$\{\varepsilon_j\}_{j\in \ZZ/3\ZZ}$ be the standard basis of the weight
lattice of $T$.
Set $\theta = 2 \varepsilon_0 - \varepsilon_1 - \varepsilon_2$, a weight
of $T$. Let $\tilde{R} = \bigoplus_{n \ge 0} R^T_{n \theta}$ be a graded
ring, where $R^T_{n \theta}$ is the semi-invariant subspace belonging to
the weight $n \theta$. Consider the projective scheme $\mathrm{Proj} \tilde{R}$
over the affine scheme $\mathrm{Spec} R^T_0$.
It is easy to see that $R^T_0$ is generated by $u_1 u_2 u_3$, $v_1 v_2 v_3$
and $u_1 v_1$ and the morphism $\CC[S] \longrightarrow R^T_0$ defined by
$a \mapsto u_1 u_2 u_3$, $b \mapsto v_1 v_2 v_3$ and $h \mapsto u_1 v_1$
is an isomorphism of algebras. Thus we have 
$\mathrm{Spec} R^T_0 \simeq S \simeq \CC^2 / (\ZZ/3\ZZ)$.
Moreover 
it is well-known that $\mathrm{Proj} \tilde{R}$ is the minimal resolution of
it, hence it is isomorphic to $\tilde{S}$.
We identify $\tilde{S}$ with $\mathrm{Proj} \tilde{R}$ and
$S$ with $\mathrm{Spec} R^T_0$

We can describe local structure of $\tilde{S}$ by the graded algebra
$\tilde{R}$. First it is easy to check that $\tilde{R}$ is generated
by the elements $u_1 u_2^2$, $v_2 v_3^2$ and $u_1 v_3 \in R^T_{\theta}$
over $R^T_0$. Consider the following open subsets of $\tilde{S}$:
\[
 U_1 = D_+(v_2 v_3^2), \qquad U_2 = D_+(u_1 v_3), \qquad 
 U_3 = D_+(u_1^2 u_2).
\]
Set
\begin{align*}
 {x}_1 = \frac{u_1 v_3}{v_2 v_3^2} ,& & {x}_2 = \frac{u_1^2 u_2}{u_1 v_3},
 & & {x}_3 = u_1 u_2 u_3, \\
 {\partial}_1 = v_1 v_2 v_3 ,& & {\partial}_2 = \frac{v_2 v_3^2}{u_1 v_3},
 & & {\partial}_3 = \frac{u_1 v_3}{u_1^2 u_2}. 
\end{align*}
Then the elements $x_i$, $\partial_i$ are local section on $U_i$.
Clearly, for $i=1$, $2$, $3$, $U_i$ is isomorphic to $\CC^2$
with coordinates $(x_i, \partial_i)$ and $\{U_i\}_{i=1,2,3}$ is an affine
open covering of $\tilde{S}$. Moreover, it is known that
$(x_i; \partial_i)$ gives a symplectic coordinate of $U_i$; i.e.
the symplectic form $\omega$ of $\tilde{S}$ coincide with 
$d x_i \wedge d \partial_i$ on $U_i$ up to scalar multiplication.
It is also clear that, for $i=1$, $2$, we have the following relations
on $U_i \cap U_{i+1}$ between the local coordinates:
\begin{equation}
\label{eq:Slod-var-local-relation} 
{x}_{i} = {\partial}_{i+1}^{-1}, \quad
{x}_{i+1} = {x}_{i}^{2} {\partial}_{i},  \quad
{\partial}_{i} = {\partial}_{i+1}^2 {x}_{i+1}.
\end{equation}
Indeed the union $U_i \cup U_{i+1}$ is isomorphic to $T^* \PP^1$ as
a symplectic manifold.

A localization of the deformed preprojective algebra of type $A$ (or
the rational Cherednik algebra of type $\ZZ / \ell \ZZ$) was constructed
in \cite{BK}. The quantized local structure with respect to the
above affine covering was studied in \cite{Kuw}.
\subsubsection{ }

Now we discuss construct an
 ACDO $\calD_{\tilde{S}, \hbar}^{ch}$ on $\tilde{S}$ as follows.
First define $\calD_{\tilde{S}, \hbar}^{ch}(U_i)$ 
to be the $\beta \gamma$-system with generators $\{x_i, \widehat{\partial}_i\}$ 
for $i=1$, $2$, $3$.
Namely $\calD_{\tilde{S}, \hbar}^{ch}(U_i)$ is an asymptotic vertex algebra which is
isomorphic to 
$\CC((\hbar))[(x_i)_{(-1)}, (\widehat{\partial}_i)_{(-1)}, (x_i)_{(-2)}, (\widehat{\partial}_i)_{(-2)}, \dots]$ 
as a vector space with a fundamental OPE
$\widehat{\partial}_i(z) x_i(w) \sim \hbar / (z-w)$. 

Next we consider a quantized chiral analogue of 
the relation (\ref{eq:Slod-var-local-relation}). For $i=1$, $2$, 
we define relations on $U_i \cap U_{i+1}$ as follows:
\begin{align}
 \label{eq:chiral-local-relation}
x_{i}(z) &= \widehat{\partial}_{i+1}(z)^{-1}, \\
x_{i+1}(z) &= \NO x_{i}(z)^2 \widehat{\partial}_{i}(z) \NO + 2 \hbar \partial_z x_{i}(z), 
 \nonumber \\
\widehat{\partial}_{i}(z) &= \NO \widehat{\partial}_{i+1}(z)^2 x_{i+1}(z) \NO 
- 2 \hbar \partial_z \widehat{\partial}_{i+1}(z).
 \nonumber 
\end{align}
Note that we need anomaly which is essentially the same one 
for $\calD^{ch}_{\PP^1}$ in the above relations. It is easy to check that
the above $\{\calD_{\tilde{S}, \hbar}^{ch}(U_i) \}_{i=1,2,3}$ and the relations
(\ref{eq:chiral-local-relation}) define a well-defined sheaf of asymptotic vertex algebras
$\calD_{\tilde{S}, \hbar}^{ch}$.

Set 
\begin{align*}
 G^+(z) &= - \hbar^{-2} \left\{ \NO x_2(z) \widehat{\partial}_2(z)^2 \NO 
 + 2 \hbar \partial_z \widehat{\partial}_2(z) \right\} 
(= - \hbar^{-2} \widehat{\partial}_1(z)), \\
 G^-(z) &= - \hbar^{-1} \left\{\NO x_2(z)^2 \widehat{\partial}_2(z) \NO 
 - 2 \hbar \partial_z \widehat{\partial}_2(z) \right\} 
(= - \hbar^{-1} x_3(z)),  \\
 J(z) &= - \hbar^{-1} \NO x_2(z) \widehat{\partial}_2(z) \NO.
\end{align*}
 They are $\rho$-invariant
global sections of $\calD_{\tilde{S}, \hbar}^{ch}$. Direct calculation shows that
OPEs of these global sections are given as follows:
\begin{gather*}
 J(z) J(w) \sim - \frac{1}{(z-w)^2}, \qquad G^{\pm}(z) G^{\pm}(w) \sim 0, \qquad
 J(z) G^{\pm}(w) \sim \pm \frac{1}{z-w} G^{\pm}(w),  \\
 G^{+}(z) G^{-}(w) \sim \frac{6}{(z-w)^3} - 
 \frac{6}{(z-w)^2} J(w) + \frac{1}{(z-w)} \{ 3 \NO J(w)^2 \NO
 - 3 \partial_w J(w) \} .
\end{gather*}
These OPEs coincide with those in (\ref{eq:ope-w23}) at critical level
 $k=- h^{\vee} = -3$. 

Note that these relations are a quantization of the vertex Poisson algebra structure (\ref{vert-poiss-on-jet-w23}). To see this one can, for example,
introduce $\hat{J}=\hbar J$, $\hat{G}^+=\hbar^2 G^+$, $\hat{G}^-=\hbar G^-$, and then we have the following relation between these fields:
\[
 \NO \hat{J}(z)^3 \NO + \NO \hat{G}^+(z) \hat{G}^-(z) \NO
 = \frac{3}{2} \hbar \partial_z \NO \hat{J}(z)^2 \NO - \hbar^2 \partial_z\hat{J}(z).
\]
This relation gives a quantization of the defining relation of the vertex 
Poisson algebra $\CC[J_{\infty} S]$
upon identification $h=\hat{J}$, $a=\hat{G}^+$, $b= - \hat{G}^-$.
Moreover,
the prescription of sect.~\ref{quasiclass-lim-hadic-vert}
will produce
relations (\ref{vert-poiss-on-jet-w23}).

\footnotesize{T.A.:  RIMS,
Kyoto University, Kyoto 606-8502 JAPAN.
Email:    arakawa@kurims.kyoto-u.ac.jp

T.K.: Faculty of Mathematics, National Research University -- Higher School Economics, 7 Vavilova Street, Moscow, 117312, Russia.
E-mail: toshiro.kuwa@gmail.com, tkuwabara@hse.ru

F.M.: Department of Mathematics, University of Southern California, Los Angeles, CA 90089, USA.
E-mail:fmalikov@usc.edu}

 \end{document}